\documentclass[11pt]{article}
\usepackage{amsxtra,amssymb,amsmath,enumerate,graphicx,bbm,hyperref,amsthm,color}
\usepackage[active]{srcltx}
\usepackage{t1enc}
\usepackage{aurical}
\usepackage[T1]{fontenc}
\usepackage{float}
\usepackage[linesnumbered,ruled,vlined]{algorithm2e}
\usepackage{dsfont}
\usepackage{url} 
\usepackage{subcaption}
\usepackage{titlesec}
\usepackage{float} 
\titleformat{\paragraph}
{\normalfont\normalsize\bfseries}{\theparagraph}{1em}{}
\titlespacing*{\paragraph}
{0pt}{3.25ex plus 1ex minus .2ex}{1.5ex plus .2ex}

\newtheorem{lemma}{Lemma}[section]
\newtheorem{proposition}[lemma]{Proposition}
\newtheorem{theorem}[lemma]{Theorem}

\newtheorem{remark}[lemma]{Remark}

\newtheorem{assumption}[lemma]{Assumption}

\newcommand{\e}{\mathbb{E}}
\newcommand{\E}{\mathbb{E}}

\renewcommand{\P}{\mathbb{P}}
\newcommand{\R}{\mathbb{R}}

\newcommand{\1}{{\mathbf 1}}

\def \Esp#1{{\mathbb E}\left[#1\right]}

\def\be{\begin{eqnarray}}
\def\ee{\end{eqnarray}}
\def\b*{\begin{eqnarray*}}
\def\e*{\end{eqnarray*}}

\def\vp{\varphi}

\def\be{\begin{eqnarray}}
\def\ee{\end{eqnarray}}
\def\beq{\begin{equation}}
\def\eeq{\end{equation}}
\def\b*{\begin{eqnarray*}}
\def\e*{\end{eqnarray*}}
\def\bi{\begin{itemize}}
\def\ei{\end{itemize}}

\def \1{{\bf 1}}
\def\vp{\varphi}
\def\eps{\varepsilon}

\def\={\;=\;}

\def\x{\times}

\def\Esp#1{\mathbb{E}\left[#1\right]}

\def \proof{{\noindent \bf Proof. }}
\def \ep{\hbox{ }\hfill$\Box$}
 \def\reff#1{{\rm(\ref{#1})}}

 \def\vs#1{\vspace{#1mm}}

\def\ti{{t_i}}
\def\tip{ {t_{i+1}} }

\def \E{\mathbb{E}}
\def \F{\mathbb{F}}

\def \M{\mathbb{M}}
\def \N{\mathbb{N}}
\def \P{\mathbb{P}}

\def \R{\mathbb{R}}
\def\T{\mathbb{T}}

\def\Bc{{\cal B}}

\def\Ec{{\cal E}}
\def\Fc{{\cal F}}

\def\Kc{{\cal K}}

\def\Pc{{\cal P}}

\def\Lb{{\mathbf L}}

\def\Yb{\bar Y}
\def\Zb{\bar Z}

\def\Esph#1{\hat \E\left[#1\right]}
\def\Xbf{\mathbf{X}}

\begin{document}

\title{Numerical approximation of  BSDEs using local polynomial drivers and branching processes}

\author{Bruno Bouchard\footnote{Universit\'e Paris-Dauphine, PSL Research University, CNRS, UMR [7534], CEREMADE, 75016 PARIS, FRANCE.}  \thanks{bouchard@ceremade.dauphine.fr} 
\and Xiaolu Tan\addtocounter{footnote}{-1}\footnotemark[\value{footnote}]~\addtocounter{footnote}{1}\thanks{tan@ceremade.dauphine.fr}
\and Xavier Warin\footnote{EDF R\&D \& FiME, Laboratoire de Finance des March\'es de l'Energie} \thanks{xavier.warin@edf.fr} 
\and Yiyi Zou\addtocounter{footnote}{-1}\footnotemark[\value{footnote}]~\addtocounter{footnote}{2}\thanks{zou@ceremade.dauphine.fr} 
}

\maketitle
 
\begin{abstract}
We propose a new numerical scheme for   Backward Stochastic Differential Equations based on branching processes. 
We approximate an arbitrary (Lipschitz) driver by local polynomials and then use a Picard iteration scheme. 
Each step of the Picard iteration can be solved by using a representation in terms of branching diffusion systems, thus avoiding the need for a fine time discretization.  
In contrast to the previous literature on the numerical resolution of BSDEs based on branching processes, we prove the convergence of our numerical scheme without limitation on the time horizon.
Numerical simulations are provided to illustrate the performance of the algorithm.
\end{abstract}

\vspace{2mm}

\section{Introduction}

	Since the seminal paper of Pardoux and Peng \cite{PardouxPeng}, the theory of  Backward Stochastic Differential Equations (BSDEs hereafter) has been largely developed,
	and has lead to  many applications in optimal control, finance, etc. (see e.g. El Karoui, Peng and Quenez \cite{ElKarouiBSDE}).
	Different approaches have been proposed during the last decade to solve them numerically, without relying on pure PDE based resolution methods. 	A first family of numerical schemes, based on a time discretization technique, has been   introduced by Bally and Pag\`es \cite{BallyPages}, Bouchard and Touzi \cite{BouchardTouzi} and Zhang \cite{Zhang}, and generated a large stream of the literature.
	The implementation of these numerical schemes requires the estimation of a sequence of conditional expectations, which can be done by  using simulations combined  with either non-linear regression techniques  or Malliavin integration by parts based representations of conditional expectations, or  by using a quantization approach, see e.g.~\cite{BouchardWarin,GLW} for references and error analysis.

	\vspace{1mm}

	Another type of  numerical algorithms is based on  a pure forward simulation of branching processes, and was  introduced by Henry-Labord\`ere \cite{Henry-Labordere_branching}, and Henry-Labord\`ere, Tan and Touzi \cite{HTT}
	(see also the recent extension by Henry-Labord\`ere et al.~\cite{HOTTW}).
	The main advantage of this new algorithm is that it avoids the estimation of   conditional expectations.
	It relies on the probabilistic representation in terms of branching processes of the 
	 so-called KPP (Kolmogorov-Petrovskii-Piskunov) equation:
	\begin{equation} \label{eq:PDEsemiLinear}
		\partial_t u(t,x) ~+~ \frac{1}{2} D^2 u(t,x) ~+~  \sum_{k \ge 0} p_k u^k(t,x)
		~=~ 0,~~~u(T,x) = g(x).
	\end{equation}
	Here, $D^2$ is the Laplacian on $\R^d$, and $(p_k)_{k\ge 0}$ is a probability mass sequence, i.e. $p_k \ge 0$ and $\sum_{k\ge 0}p_k = 1$.
	This is a  natural extension of the classical Feynmann-Kac formula, which is  well known 
	since the works of Skorokhod \cite{Skorokhod}, Watanabe \cite{Watanabe} and McKean \cite{McKean_1975}, among others. The PDE \reff{eq:PDEsemiLinear} corresponds to a BSDE with a polynomial driver and terminal condition $g(W_{T})$:
	$$
	Y_{\cdot} =g(W_{T})+\int_{\cdot}^{T} \sum_{k \ge 0} p_k (Y_{t})^{k}dt -\int_{\cdot}^{T}Z_{t}dW_{t},
	$$
	 in which $W$ is a Brownian motion. Since $Y_{\cdot} = u(\cdot,W_{\cdot})$, the $Y$-component of this BSDE can be estimated by making profit of the branching process based Feynman-Kac representation of \reff{eq:PDEsemiLinear} by means of a pure forward Monte-Carlo scheme, see Section \ref{sec: repre} below. The idea is not new. It was already proposed in Rasulov, Raimov and Mascagni \cite{RasulovRaimovMascagni}, although no rigorous convergence analysis was provided.	Extensions to more general drivers can be found in \cite{Henry-Labordere_branching,HOTTW,HTT}.  Similar algorithms have been studied by Bossy et al.~\cite{BCLMVY}  to solve non-linear Poisson-Boltzmann equations. 
	
	\vspace{1mm}
	It would be tempting to use this representation to solve BSDEs with Lipschitz drivers, by approximating their drivers by  polynomials. This is however not  feasible in general. The  reason is that PDEs (or BSDEs) with polynomial drivers, of degree bigger or equal to two, typically explode in finite time. They are only well posed on a small time interval. It is worse when the degree of the polynomial increases. Hence, no convergence can be expected for the case of general drivers.

	\vspace{1mm}

	In this paper, we propose to instead use a local polynomial approximation.  Then, convergence of the sequence of approximating drivers to the original one can be ensured without explosion of the corresponding BSDEs, that can be defined on a arbitrary time interval. It requires to be combined with a Picard iteration scheme, as the choice of the  polynomial form will depend on the position in space of the solution $Y$ itself. However, unlike classical Picard iteration schemes for BSDEs, see e.g.~Bender and Denk \cite{BenderDenk}, we do not need to have a very precise estimation of the whole path of the solution at  each Picard iteration. Indeed, if local polynomials are fixed on a partition $(A_{i})_{i}$ of $\R$, then one only needs to know in which $A_{i}$ the solution stays at certain branching times of the underlying branching process. If the $A_{i}$'s are large enough, this does not require a very good precision in the intermediate estimations. We refer to Remark \ref{rem: localisation} for more details. 

	\vspace{1mm}
	
	We finally insist on the fact that our results will be presented in a Markovian context for simplification. However, all of our arguments work trivially in a non-Markovian setting too.

\section{Numerical method for a class of BSDE based on branching processes}
\label{sec:model}

	Let $T>0$, $W$ be a standard $d$-dimensional Brownian motion on a filtered probability space $(\Omega, \Fc, \F = (\Fc_t)_{t\ge 0}, \P)$,
	and $X$ be the solution of the stochastic differential equation:
	\begin{equation} \label{eq: Diffusion}
		{X}=X_{0}+\int_{0}^{\cdot} \mu({X}_s)\,dt+\int_{0}^{\cdot} \sigma( {X}_s)\,dW_s,
	\end{equation}
	where $X_0$ is a constant, lying  in a compact subset $\Xbf$ of $\R^{d}$, and $(\mu,\sigma):[0,T]\x \R^{d}\mapsto \R^{d}\x \M^{d}$ is assumed to be Lipschitz continuous with  support contained in $\Xbf$. 
	Our aim is to provide a numerical scheme for the resolution of the backward stochastic differential equation 
	\be\label{eq: BSDE}
		Y_{\cdot}=g({X}_T)+\int_{\cdot}^{T}f( {X}_s,Y_s)\,ds-\int_{\cdot}^{T}Z_s\,dW_s.
	\ee
	In the above, $g:\R^{d}\mapsto \R$ is assumed to be  measurable and  bounded, $f\in \R^{d}\x \R\mapsto \R$ is measurable with linear growth and Lipschitz in its second argument, uniformly in the first one. As a consequence, there exists $M\ge 1$ such that 
 \be\label{eq: borne X Y g}
 |g(X_{T})|\le M \;\mbox{ and } |X|+|Y|\le M\;\mbox{ on } [0,T].
 \ee 
 
 \begin{remark} The above conditions are imposed to easily localize the solution $Y$ of the BSDE, which will be used in our estimates later on. One could also assume that $g$ and $f$ have polynomial growth in their first component and that $\Xbf$ is not compact. After possibly truncating the coefficients and reducing their support, one would go back to our conditions. Then, standard estimates and stability results for SDEs and BSDEs could be used to estimate the additional error in a very standard way.  See e.g.~\cite{ElKarouiBSDE}.
 \end{remark}

\subsection{Local polynomial approximation of the generator}

	A first main ingredient of our algorithm consists in approximating the driver $f$ by a driver $f_{\ell_{\circ}}$ 
	that has a local polynomial structure. 
	Namely, let 
	\be\label{eq: f loc poly}
		f_{\ell_{\circ}}: (x,y,y')\in   \R^{d}\x \R\x \R\mapsto \sum_{j=1}^{j_{\circ}} \sum_{\ell=0}^{\ell_\circ}a_{j,\ell}(x)y^{\ell}\vp_{j}(y'),
	\ee 
in which $(a_{j,\ell},\vp_{j})_{\ell \le \ell_{\circ}, j\le j _{\circ}}$ is a family of continuous and bounded maps satisfying 
	\be\label{eq: borne a vp}
		|a_{j,\ell}  | \le C_{\ell_{\circ}} \;,\; |\vp_{j}(y'_{1})-\vp_{j}(y'_{2})|\le L_{\varphi}|y'_{1}-y'_{2}| \mbox{ and } |\vp_{j} |\le 1, 
	\ee 
for all $y'_{1},y'_{2}\in \R$, $j\le j_{\circ}$ and $\ell\le \ell_{\circ}$,  for some constants $C_{\ell_{\circ}},L_{\varphi}\ge 0$. In the following, we shall assume that $\ell_{\circ}\ge 2$ (without loss of generality).  One can think of the $(a_{j,\ell})_{\ell \le \ell_{\circ}}$ as the coefficients of a polynomial 
approximation of $f$ on a subset $A_{j}$, the $A_{j}$'s forming a partition of   $[-M,M]$. Then, the $\vp_{j}$'s have to be considered as smoothing kernels that allow one to pass in a Lipschitz way from one part of the partition to another one. We therefore assume that 
\be\label{eq: intersection vp}
	{  \#\{j\in \{1,\cdots,j_{\circ}\}: \; \vp_{j}(y)>0\}\le 2 \;\mbox{ for all } y \in \R, }
\ee
and that  $y \mapsto f_{\ell_{\circ}}(x, y, y)$ is globally Lipschitz. In particular, 
\be\label{eq: BSDElo}
	\Yb_{\cdot}
	=
	g({X}_T)
	+
	\int_{\cdot}^{T}f_{\ell_{\circ}}( {X}_s,\Yb_s,\Yb_s)\,ds
	-\int_{\cdot}^{T} \Zb_s\,dW_s,
\ee
has a unique solution $( \Yb,  \Zb)$ such that $\E[\sup_{[0,T]} |\bar Y|^{2}]<\infty$.
Moreover, by standard estimates, $(\Yb, \Zb)$ provides a good approximation of $(Y,Z)$ whenever $f_{\ell_{\circ}}$ is a good approximation of $f$:
\be \label{eq: erreur driver}
	\E \Big[\sup_{[0,T]}|Y- \Yb |^{2} \Big]+\E\Big[\!\! \int_{0}^{T} \!\!\!\!|Z_{t}- \Zb_{t}|^{2} dt \Big]
	\le
	C \E\Big[\!\! \int_{0}^{T} \!\!\! |f-f_{\ell_{\circ}}|^{2}(X_{t},Y_{t},Y_{t})dt\Big],
\ee
for some $C>0$ that does not depend on $f_{\ell_{\circ}}$, see e.g.~\cite{ElKarouiBSDE}.

The choice of $f_{\ell_{\circ}}$ will obviously depend on the application at hand and does not need to be more commented. Let us just mention that our algorithm will be more efficient if the sets $\{y \in \R : \vp_{j}(y) = 1\}$  are large and the intersection between the supports of the $\vp_{j}$'s are small, see Remark \ref{rem: localisation} below.  

{We also assume that 
\begin{align}\label{eq: bone bar Y}
|\bar Y| \le M.
\end{align}
Since we intend to keep $f_{\ell_{\circ}}$ with linear growth in its first component, and bounded in the two other ones, uniformly in $\ell_{\circ}$, this is without loss of generality. 
}
\subsection{Picard iteration with doubly reflected BSDEs}

	Our next step is to introduce a Picard iteration scheme to approximate the solution $\Yb$ of  \eqref{eq: BSDElo}.
	Note however that, although  the map $y \mapsto f(x, y, y)$ is globally Lipschitz,   the map $y \mapsto f(x, y, y')$ is a polynomial, given $y'$, and hence only locally Lipschitz in general.
	In order to reduce to a Lipschitz driver, we shall  apply our Picard scheme to  a  doubly (discretely) reflected BSDE, with   lower and upper barrier given  by the bounds $-M$ and $M$ for $\Yb$, recall \reff{eq: bone bar Y}. 

	\vs2

	Let $h_{\circ}$ be defined by \eqref{eq:def_delta} in the Appendix. It is a lower bound for the explosion time of the BSDE with driver $y \mapsto f(x, y, y')$. Let us then fix $h \in (0, h_{\circ})$ such that $N_{h}:=T/h \in \N$, and define  
	\begin{align}\label{eq: def ti}
		\ti=ih
		\;\;\;\;\;\mbox{and}\;\;\;\;
		\T := \{\ti, ~i= 0, \cdots, N_{h}\}.
	\end{align}
	We initialize our  Picard scheme by setting 
	\be\label{eq: borne Yn0}
		\Yb^{0}_t ={\rm y}(t, X_t)\;\mbox{ for } t \in [0,T],
	\ee 
	in which ${\rm y}$ is a deterministic function, bounded by $M$ and such that ${\rm y}(T,\cdot)=g$.
	Then, given  $\Yb^{m-1}$, for $m\ge 1$, we define $(\Yb^m, \Zb^m, \bar K^{m,+}, \bar K^{m,-})$ as the solution on $[0,T]$ of 
	\be
		&&\Yb^{m}_t
		=
		g(X_T) +\int_t^{T}{f}_{\ell_{\circ}}(X_s,\Yb^{ m}_s, \Yb^{m-1}_s)\,ds
		-\int_t^{T} \Zb^{m}_s\,dW_s
		+ \int_{[t,T] \cap \T} \!\!\!\!\! d (\bar K^{m,+} -\bar K^{m,-})_s,\nonumber
		\;\;\;\\
		&&-M \le \Yb^m_t \le M, ~~~\forall t \in \T,~a.s. \label{eq: IterationBSDE}\\
		&&\int_{\T} (\Yb^m_s + M) d\bar K^{m,+}_s = \int_{\T} (\Yb^m_s - M) d\bar K^{m,-}_s = 0,\nonumber
	\ee
	where $\bar K^{m,+}$ and $\bar K^{m,-}$ are non-decreasing processes.

	\begin{remark}
		  Since the solution $\Yb$ of   \eqref{eq: BSDElo} is bounded by $M$, 
		the quadruple of processes $(\Yb,  \Zb, \bar K^+, \bar K^-)$ (with $\bar K^+ \equiv \bar K^- \equiv 0$) is in fact the unique solution of the same reflected BSDE as in \eqref{eq: IterationBSDE} but 
		with   $f_{\ell_{\circ}}(X, \Yb, \Yb)$ in place of ${f}_{\ell_{\circ}}(X_s,\Yb^{ m}, \Yb^{m-1})$.

	\end{remark}

	\begin{remark} \label{rem:DRBSDE_discret}
		 One can equivalently define  the process $\Yb^m$ in a recursive way.
		Let $\Yb^m_T := g(X_T)$ be the terminal condition, and define, on each interval $[t_i, t_{i+1}]$,   $(Y^m_{\cdot}, Z^m_{\cdot})$ as the solution on $[t_{i},t_{i+1}]$ of
		\be \label{eq:def_Ym}
			Y^m_{\cdot} 
			~=~
			\Yb^m_{t_{i+1}}
			 +
			 \int^{t_{i+1}}_{\cdot} \!\!\!\! f_{\ell_{\circ}} (X_s, Y^m_s, \Yb^{m-1}_s) ds - \int^{t_{i+1}}_{\cdot} Z^m_s dW_s.\;\;
		\ee
		Then,  $\Yb^m  := Y^m $ on $ (t_i, t_{i+1}]$, and $\Yb^m_{t_i} := (-M) \vee Y^m_{t_i} \wedge M$.
	\end{remark}
	
	The error due to our Picard iteration scheme is handled in a standard way. It depends on the constants
	\b*
		L_1 := 2C_{\ell_{\circ}} \sum_{\ell =1}^{\ell_{\circ}} \ell (M_{h_{\circ}})^{\ell -1} ,
		\;\;\;M_{h_\circ}
		 L_2 := L_{\varphi} \sum_{\ell=0}^{\ell_{\circ}} 2 C_{\ell_{\circ}} (M_{h_{\circ}})^{\ell},
	\e*
	where $M_{h_\circ}$ is defined by \eqref{eq:def_M_delta}.
	
	\begin{theorem}\label{thm: main} 
	The system \reff{eq: IterationBSDE} admits a unique solution $(\bar Y^{m},\bar Z^{m}, \bar K^{m,+}, \bar K^{m,-})_{m\ge 0}$ such that 
	$\bar Y^m$ is uniformly bounded for each $m\ge 0$. Moreover, for all   $m\ge 0$, 
		 $|\bar Y^{m}|$ is   uniformly bounded by the constant $M_{h_{\circ}}$, and 
		\begin{align*}
			 |\Yb^{m}_t- \Yb_t |^2
			~\le~
			\frac{L_2 }{\lambda^2}\Big(\frac{L_2 (T-t)}{\lambda^2} \Big)^m (2M)^2 \frac{e^{\beta T}}{\beta},
		\end{align*}
		 for all $t\le T$, and all constants $\lambda >0$, $\beta > 2L_1 + L_2 \lambda^2$.
	\end{theorem}
	\proof
	$\mathrm{i)}$ 
	First, when $\Yb^m$ is uniformly bounded, $f_{\ell_{\circ}}(X_s, \Yb^m_s, \Yb^{m-1}_s)$ can be considered to be uniformly Lipschitz in $\Yb^{m}$,
	then  \reff{eq: IterationBSDE}  has at most one bounded solution.
	Next, in view of Lemma \ref{lemm:estim_ODE} and Remark \ref{rem:DRBSDE_discret}, 
	it is easy to see that \eqref{eq:def_Ym} has a unique solution $Y^m$, bounded by $M_{h_{\circ}}$ (defined by \eqref{eq:def_M_delta}) on each interval $[t_i, t_{i+1}]$.
	It follows the existence of the solution to \eqref{eq: IterationBSDE}.
	Moreover, $\Yb^m$ is also bounded by $M_{h_{\circ}}$ on $[0,T]$, and more precisely bounded by $M$ on the discrete grid $\T$, by construction.

	\vspace{1mm}
	
	$\mathrm{ii)}$ 	Consequently, the generator $f_{\ell_{\circ}}(x,y, y')$ can be considered to be uniformly Lipschitz in $y$ and $y'$.
	Moreover, using \eqref{eq: borne a vp} and \eqref{eq: intersection vp}, one can identify the corresponding Lipschitz constants as $L_1$ and $L_2$.	
	
	\vspace{1mm}

	Let us denote $\Delta \Yb^m := \Yb^m - \Yb$ for all $m \ge 1$.
	We notice that, in Remark \ref{rem:DRBSDE_discret}, the truncation operation $\Yb^m_{t_i} := (-M) \vee Y^m_{t_i} \wedge M$ can only make the value $(\Delta \Yb^m_{{t_{i}}})^2$ smaller than $(Y^m_{t_i} - \Yb_{t_i})^2$, since $|\bar Y|\le M$.
	Thus we can apply It\^o's formula to $(e^{\beta t} (\Delta \Yb^{m+1}_{t})^2)_{t\ge 0}$ on each interval $[t_i, t_{i+1}]$,
	and then take expectation to obtain
	\b*
		&&
		\E \big[ e^{\beta t} (\Delta \Yb^{m+1}_t)^2 \big]
		+ \beta \E \Big[ \int_t^T e^{\beta s} |\Delta \Yb^{m+1}_s|^2 ds 
		+ \int_t^T e^{\beta s} | \Delta \Zb^{m+1}_s|^2 ds \Big] \\
		&\le&
		2 \E \Big[ \int_t^T e^{\beta s} \Delta \Yb^{m+1}_s \big( f_{\ell_{\circ}} (X_s, \Yb^{m+1}_s, \Yb^m_s) - f_{\ell_{\circ}}(X_s, \Yb_s, \Yb_s) \big) ds \Big].
	\e*
	Using the Lipschitz property of $f_{\ell_{\circ}}$ and the inequality $\lambda^2 + \frac{1}{\lambda^2} \ge 2$,
	it follows that the r.h.s.~of the above inequality is bounded by
	$$
		(2L_1 + L_2 \lambda^2) \E \Big[ \int_t^T e^{\beta s} (\Delta \Yb^{m+1}_s)^2 ds \Big]
		+
		\frac{L_2}{\lambda^2} \E \Big[ \int_t^T e^{\beta s} (\Delta \Yb^m_s)^2 ds \Big].
	$$
	Since $\beta \ge 2L_1 + L_2 \lambda^2$,  the above implies 
	\be \label{eq:bound_DeltaY_t}
		\E \Big[ e^{\beta t} (\Delta \Yb_t^{m+1})^2 \Big] 
		~\le~
		\frac{L_2}{\lambda^2} \E \Big[\int_t^T e^{\beta s} (\Delta \Yb^m_s)^2 ds \Big],
	\ee
	and hence
	\b*
		\E \Big[ \int_0^T e^{\beta t} (\Delta \Yb_t^{m+1})^2 dt \Big] 
		&\le&
		\frac{L_2}{\lambda^2} T \E \Big[\int_0^T e^{\beta s} (\Delta \Yb^m_s)^2 ds \Big] .
	\e*
	Since $|\Delta \Yb^{0}|=|{\rm y}(\cdot, X)-\bar Y|\le 2M$ by \reff{eq: bone bar Y} and our assumption $|{\rm y}|\le M$, this shows that 
		\b*
		\E \Big[ \int_0^T e^{\beta t} (\Delta \Yb_t^{m})^2 dt \Big] 
		&\le&
		\big( \frac{L_2}{\lambda^2}T\big)^{m} (2M)^2 e^{\beta T}/\beta.
	\e*
	Plugging this in \eqref{eq:bound_DeltaY_t} leads to the required result at $t=0$. 
	It is then clear that the above estimation does not depend on  the initial condition $(0, X_0)$, so that the same result holds  true for every $t \in [0,T]$.
	\qed

\subsection{A branching diffusion representation for $\Yb^m$}\label{sec: repre}

	We now explain how the  solution of   \eqref{eq:def_Ym} on $[t_i, t_{i+1})$ can be represented by means of a branching diffusion system.
	More precisely, let  us  consider an element  {$(p_\ell)_{0\le \ell\le \ell_{\circ}} \in \R^{\ell_{\circ}+1}_+$ such that $\sum_{\ell = 0}^{\ell_{\circ}} p_{\ell} = 1$}, set $K_{n}:=\{(1,k_{2},\ldots,k_{n}): (k_{2},\ldots,k_{n})\in \{0,\ldots,\ell_{\circ}\}^{n}\}$ for $n\ge 1$, and $K:=\cup_{n\ge 1} K_{n}$. Let $(W^{k})_{k\in K }$ be a sequence of independent $d$-dimensional Brownian motions, $(\xi_{k})_{k\in  K }$ and $(\delta_{k})_{k\in  K}$ be two sequences of independent random variables, such that 
$$
\P[\xi_{k}=\ell]=p_{\ell}, \,\;\ell \le \ell_{\circ}, k\in  K , 
$$
and 
\be\label{eq: def bar F} 
\bar F(t):=\P[\delta_{k}>t]=\int_{t}^{\infty} \rho(s)ds,\; t\ge 0,\;k\in K , 
\ee
for some continuous strictly positive map $\rho: \R_+ \to \R_+$. We assume that 
\be\label{hyp eq: independance} 
\mbox{$(W^{k})_{k\in K }$, $(\xi_{k})_{k\in K }$, $(\delta_{k})_{k\in K }$ and $W$ are independent.}
\ee
Given the above, we construct particles $X^{(k)}$ that have the dynamics \reff{eq: Diffusion} up to a killing time $T_{k}$ at which they split in $\xi_{k}$ different (conditionally) independent particles with dynamics   \reff{eq: Diffusion} up to their own killing time. The construction is done as follows. First, we set 
  $T_{(1)}:=\delta_{1}$, and, given $k=(1,k_{2},\ldots,k_{n}) \in K_{n}$ with $n\ge 2$, we let $T_{k}:=\delta_{k}+T_{k-}$ in which $k-:=(1,k_{2},\ldots,k_{n-1}) \in K_{n-1}$. 
  By convention, $T_{(1)-}=0$.  We can then define the Brownian particles $(W^{(k)})_{k\in K}$ by using the following induction. We first set 
$$
W^{((1))}:=W^{1}\1_{[0,T_{(1)}]}\;,\;\Kc^{1}_{t}:=\{(1)\}\1_{[0,T_{(1)})}(t)+\emptyset  \1_{[0,T_{(1)})^{c}}(t),
~~\mbox{and}~\bar \Kc^1_t =\{(1)\} \;t\ge 0.
$$
Then, given $n\ge 2$, we define 
\begin{equation*}
W^{(k\oplus j)}:=\left(W^{(k)}_{\cdot\wedge T_{k}} +(W^{{k\oplus j}}_{\cdot \vee  T_{k}}-W^{{k\oplus j}}_{T_{k}})\1_{j\ne 0}\right)\1_{[0,T_{k\oplus j}]}, \;  0\le j\le \xi_{k}, \; k\in K_{n-1},
\end{equation*}
and
$$
	\bar \Kc^n_{t}:=\{k\oplus j: k\in \bar \Kc^{n-1}_{T}, 1\le j\le \xi_{k} ~\mbox{s.t.}~ t\ge T_{k}] \},
	\;\;
	\bar \Kc_{t} := \cup_{n \ge 1} \bar \Kc^n_{t},
$$
\b*
	\Kc^n_t := \{ k\oplus j: k\in \bar \Kc^{n-1}_{T}, 1\le j\le \xi_{k} ~\mbox{s.t.}~ t\in [T_{k},T_{k\oplus j})\},
	\;\;
	\Kc_t := \cup_{n\ge 1} \Kc^n_t,
\e*
in which we use the notation $(1,k_{1},\ldots,k_{n-1})\oplus j= (1,k_{1},\ldots,k_{n-1},j)$.
In other words, $\bar \Kc^n_{t}$ is the collection of particles of the $n$-th generation that are born before   time $t$ on, while   $  \Kc^n_{t}$ is the collection of particles in $\bar \Kc^n_{t}$ that are still alive at $t$.

Now observe that the solution $X^{x}$ of \reff{eq: Diffusion} on $[0,T]$ with initial condition $X^{x}_{0}=x\in \R^{d}$ can be identified in law on the canonical space as a process of the form $\Phi[x](\cdot,W)$ in which the deterministic map $(x,s,\omega)\mapsto \Phi[x](s,\omega)$ is $\Bc(\R^{d}) \otimes \Pc$-measurable, where $\Pc$ is the predictable $\sigma$-filed on $ [0,T]\x \Omega$. We then define the corresponding particles $(X^{x,(k)})_{k\in K}$ by $X^{x,(k)}:=\Phi[x](\cdot,W^{(k)})$.   

\vspace{2mm}

Given the above construction, we can now    
introduce a sequence of deterministic map associated to $(\bar Y^{m})_{m\ge 0}$. First, we set 
\be\label{eq: def v0m vn0}
 v^{0}:={\rm y}\ ,
\ee
recall \reff{eq: borne Yn0}. 
Then, given $v^{m-1}$ and $v^m(t_{i+1}, \cdot)$, we define
\b*
V^{m}_{t,x}&:=& \Big(\prod_{k \in   \Kc_{\tip-t}}G^{m}_{t,x}(k) \Big) \Big(\prod_{k \in\bar{\Kc}_{\tip-t}\setminus\Kc_{\tip-t}}A^{m}_{t,x}(k) \Big),   \\
G^{m}_{t,x}(k)&:=& \frac{v^{m} \big(t_{i+1}, X_{\tip-t}^{x,(k)} \big)}{\bar{F}({\tip-t}-T_{k-})},\\
A^{m}_{t,x}(k)&:=& \frac{\sum_{j=1}^{j_{\circ}}a_{j,\xi_k}(X_{T_k}^{x,(k)})\vp_{j}(v^{m-1}(t+T_{k},X_{T_k}^{x,(k)}))}{p_{\xi_k}\,\rho(\delta_k)},
\;\;\;
\e*
$\forall (t,x)\in [\ti,\tip)\x \Xbf.$
We finally set, whenever $V^m_{t,x}$ is integrable,
\begin{equation*}
	v^{m}(t,x):=\E\left[V^{m}_{t,x}\right], \;\;(t,x)\in (t_i, t_{i+1}) \x \Xbf, ~m\ge 1,
\end{equation*}
and
\begin{equation}\label{eq: def vnm}
	v^{m}(t_i,x):= ( -M) \vee \E\left[V^{m}_{t_i,x}\right] \wedge M, \;\; x \in \Xbf, ~m\ge 1.
\end{equation}

\begin{proposition}\label{propeq: representation Ynm} 
For all $m \ge 1$ and  $(t,x) \in [0,T] \x \mathbf{X}$, the random variable $V^m_{t,x}$ is integrable.
Moreover, one has $\Yb^{m}_{\cdot} = v^{m}(\cdot,X)$ on $[0,T]$.
\end{proposition}
 
This follows from Proposition \ref{prop: representation Ynm} proved in the Appendix, which is in spirit of \cite{HOTTW}. 
The main use of this representation result here is that it provides a numerical scheme for the approximation of the component $\Yb$  of \reff{eq: BSDElo}, as explained in the next section.

\subsection{The numerical algorithm}
\label{subsec:NumAlgo}

	The representation result in Proposition \ref{propeq: representation Ynm} suggests to use a simple Monte-Carlo estimation of the expectation in the definition of $v^{m}$ based on the simulation of the corresponding particle system. However, it requires the knowledge of $  v^{m-1}$ in the Picard scheme which is used to localize our approximating polynomials. We therefore need to approximate the corresponding (conditional) expectations at each step of the Picard iteration scheme. In practice, we shall replace the expectation operator $\E$ in the definition of $v^{m}$ by an operator $\hat \E$ that can be computed explicitly, see Remark \ref{rem: cond expe} below.
	
	\vs2
	
	In order to perform a general (abstract) analysis, let us first recall  that we have defined $v^m(t,x) =$ $\E[ V_{t,x}(v^m(t_{i+1}, \cdot), v^{m-1}(\cdot) ]$ for all $t \in (t_i, t_{i+1})$
	and $v^m(t_i,x) = (-M) \vee \E[ V_{t_i,x}(v^m(t_{i+1}, \cdot), v^{m-1}(\cdot) ] \wedge M$,
	where, given two functions $\phi,\phi' : (t_i, t_{i+1}] \x \R^d \to \R$,
	\begin{align}
         \label{eq: practReg}
		V_{t,x}(\phi,\phi')&:=
		\Big(\prod_{k \in   \Kc_{\tip-t}}G_{t,x}(\phi,k) \Big) 
		\Big(\prod_{k \in\bar{\Kc}_{\tip-t}\setminus\Kc_{\tip-t}}A_{t,x}(\phi',k) \Big),\nonumber \\
		 G_{t,x}(\phi,k)&:=\frac{\phi(\tip,X_{\tip-t}^{x,(k)})}{\bar{F}({\tip-t}-T_{k-})}, \nonumber\\
		 A_{t,x}(\phi',k)&:= \frac{\sum_{j=1}^{j_{\circ}}a_{j,\xi_k}(X_{T_k}^{x,(k)})\vp_{j}(\phi'(t+T_{k},X_{T_k}^{x,(k)})}{p_{\xi_k}\,\rho(\delta_k)}.
	\end{align}	
	
	Let us then denote by $\Lb^{\infty}_{M_{h_{\circ}}}$ the class of all Borel measurable functions $\phi: [0,T] \x \R^d \to \R$ that are bounded by $M_{h_{\circ}}$,
	and let $\Lb^{\infty}_{M_{h_{\circ}}, 0} \subset \Lb^{\infty}_{M_{h_{\circ}}}$ be a subspace, generated by a finite number of basis functions.  
	Besides,  let us consider a sequence $(U_i)_{i \ge 1}$ of i.i.d.~random variables of uniform distribution on $[0,1]$, independent of $(W^k)_{k \in K}$, $(\xi_k)_{k \in K}$, $(\delta_k)_{k \in K}$ and $W$ introduced in \eqref{hyp eq: independance}.
	Denote   $\hat \Fc := \sigma(U_i, i \ge 1)$.

From now on, we use the notations
$$
\|\phi\|_{\ti}:=\sup_{(t,x)\in [\ti,\tip)\x \R^d}|\phi(t,x)| \mbox{ and }  \|\phi\|_{\infty}:=\sup_{(t,x)\in [0,T]\x \R^d}|\phi(t,x)|
$$
for all functions $\phi: [0,T] \x \R^d \to \R$.

\begin{assumption} \label{assum:hatE}
	There  exists an operator  $\hat \E[ \hat V_{t,x}(\phi, \phi')] (\omega)$, defined for all $\phi, \phi' \in \Lb^{\infty}_{M_{h_{\circ}}, 0}$,
	such that $(t,x, \omega) \mapsto \hat \E[ \hat V_{t,x}(\phi, \phi')] (\omega)$ is  $ \Bc([0,T] \x \R^d)\otimes \hat \Fc $-measurable,
	and such that the function $(t,x) \in [0,T] \x \R^d\mapsto \hat \E[ \hat V_{t,x}(\phi, \phi')](\omega)$ belongs to $\Lb^{\infty}_{M_{h_{\circ}}, 0}$
	for every fixed $\omega \in \Omega$.
	Moreover, one has
	$$
		\Ec(\hat \E)
		:=
		\|\sup_{\phi,\phi'\in \Lb^{\infty}_{M_{h_{\circ}},0}}\Esp{|\Esp{V_{\cdot}(\phi,\phi')}-\hat \E[ {\hat V_{\cdot}(\phi,\phi')}] |}\|_{\infty}
		< \infty.
	$$
\end{assumption}

	In practice, the operator $\hat \E[\hat V]$ will be decomposed in two terms: $\hat V$ is an approximation of the operator $V$ defined with respect to a finite time grid that projects the arguments $\phi$ and $\phi'$ on a finite  functional space, while $\hat \E[\hat V_{\cdot}(\cdot)]$ is a Monte Carlo estimation of $\E[\hat V_{\cdot}(\cdot)]$. See Remark \ref{rem: cond expe}.
	\vs2
	
	Then, one can construct a numerical algorithm by first setting $\hat v^{0} \equiv {\rm y}$, $\hat v^{m}(T,\cdot)=g$, $m\ge 1$, and then by defining  by induction over $m\ge 1$

	\b*
		\hat v^{m}(t,x):=(-M_{h_\circ}) \vee \Esph{\hat V_{t,x}(\hat v^{m}(\tip,\cdot),\hat v^{m-1})}\wedge M_{h_\circ}, \; t\in (\ti, \tip),
	\e*
	and
	\begin{equation}\label{eq: def hat vnm}
		\hat v^{m}(t_i,x):=(-M)\vee \Esph{\hat V_{t_i,x}(\hat v^{m}(\tip,\cdot),\hat v^{m-1})}\wedge M.
	\end{equation}

	In order to analyze the error due to the approximation of the expectation error, let us set 
	$$
		\bar q_{t}:=  \# \bar{\Kc}_{t}
		\;\;\mbox{,}\;\;
		q_{t}:=  \# \Kc_{t},
	$$
	and denote 
	$$
		V_t^M 
		:= 
		\Big( \prod_{k \in \Kc_{ {t}}} \frac{M}{\bar F( {t} - T_{k-})} \Big)
		\Big( \prod_{k \in \bar \Kc_t \setminus \Kc_t} \frac{ 2C_{\ell_{\circ}}}{p_{\xi_k} \rho(\delta_k)} \Big).
	$$
	Recall that $h<h_{\circ}$ that is defined by \eqref{eq:def_delta} in the Appendix. 
\begin{lemma}
	The two constants   
	 $$
		M^1_{h}
		:= 
		\sup_{0 \le t \le h}\E \big[  q_t V_t^M \big]
		\;\;\;\mbox{and}\;\;
		M^2_{h}
		:= 
		\sup_{0 \le t \le h}\E \big[   \bar q_t V_t^M \big]
	$$
	are finite. 
\end{lemma}
\proof
	Notice that for any constant $\eps > 0$, 
	there is some constant $C_{\eps} > 0$ such that $n \le C_{\eps} (1+ \eps)^n$ for all $n \ge 1$.
	Then 
	$$
		M^1_{h} \le C_{\eps} \E \big[ \sup_{0 \le t \le h} (1+\eps)^{q_t} V^M_t \big]
		\le C_{\eps} \E\Big[
			\prod_{k \in \Kc_h} \frac{(1+\eps)M}{\bar F(h - T_{k-})} 
			\prod_{k \in \bar \Kc_h \setminus \Kc_h} \frac{ 2 (1+\eps) C_{\ell_{\circ}}}{p_{\xi_k} \rho(\delta_k)}
		\Big],
	$$
	where the latter expectation is finite for $\eps$ small enough. This follows from the fact that  	$h < h_{\circ}$ for $h_{\circ}$ defined by \eqref{eq:def_delta} and from the same arguments as in Lemma \ref{lemm:estim_ODE} in the Appendix.
	One can similarly obtain that $M^2_{h}$ is also finite.
\qed

\begin{proposition}\label{prop: erreur discretisation et esperance} 
	Let  Assumption \ref{assum:hatE} hold true.
	Then  
	{
	\begin{align*} 
		\| \Esp{| v^{m} -\hat v^{m}|} \|_{\infty}
		&\le
		\Ec(\hat \E) 
		\big( 1+N_{h} \big)
		\frac{(m+N_{h})^{N_{h}}}{N_{h} !} 
		\Big((2L_{\varphi} M^2_{h}) \vee \frac{M^1_{h}}{M} \vee 1 \Big)^{m+N_{h}}.
	\end{align*}
	}
\end{proposition}

Before turning to the proof of the above, let us comment on the use of this numerical scheme. 
\begin{remark}\label{rem: cond expe}
In practice, the approximation of the expectation operator can  be simply constructed by using pure forward simulations of the branching process. Let us explain  this first in the case $h_{\circ}=T$. Given that $\hat v^{m}$ has already been computed, one takes it as a given function, 
one draws some independent copies of the branching process (independently of $\hat v^{m}$) and computes
$\hat v^{m+1}(t,x)$ as the Monte-Carlo counterpart of $\E[V_{t,x}(\hat v^{m+1}(T,\cdot), \hat v^{m})]$, and truncates it with the a-priori bound $M_{h_{\circ}}$ for $(\Yb^{m})_{m\ge 1}$. This corresponds to the operator $\hat \E[\hat V_{t,x}(\hat v^{m+1}(T,\cdot), \hat v^{m})]$.  If $h_{\circ}<T$, one needs to iterate backward over the periods $[t_{i},t_{i+1}]$. Obviously one cannot in practice compute the whole map $(t,x)\mapsto \hat v^{m+1}(t,x)$ and this requires an additional discretization on a suitable time-space grid. Then, the additional error analysis can be handled for instance  by using the continuity property of $v^m$ in Proposition \ref{prop: holder vm} in the Appendix. This is in particular the case if one just computes $\hat v^{m+1}(t,x)$ by replacing ${(t,x)}$ by its projection on a discrete time-space grid.
\end{remark}

\begin{remark}\label{rem: localisation} 

$\mathrm{i)}$.~In the classical time discretization schemes of BSDEs, such as those in \cite{BouchardTouzi,GLW,Zhang},
one needs to let the time step go to $0$ to reduce the discretization error.
Here, the representation formula in Proposition \ref{propeq: representation Ynm} has no discretization error related to the BSDE itself (assuming the solution of the previous Picard iteration is known perfectly), we only need to use a fixed discrete time grid $(t_i)_{0 \le i \le N_h}$ for $t_i = i h$ with $h$ small enough.

\vspace{1mm}

$\mathrm{ii)}$.~Let  $A'_{j}:=\{y \in \R ~: \vp_{j}(y) =1\}\subset A_{j}$ for $j\le j_{\circ}$, and assume that the $A'_{j}$'s are {disjoint}. If the $A'_{j}$ are large enough, we do not need to be  very precise on $\hat v^{m}$ to obtain a good approximation of $\E[V_{t,x}(  g ,   v^{m})]$ by $\E[V_{t,x}(g, \hat v^{m})]$ for $t\in [t_{N_{h}-1},t_{N_{h}})$. One just needs to ensure that $\hat v^{m}$ and $v^{m}$ belong to the same set $A'_{j}$ at the different branching times and at the corresponding $X$-positions. 
We can therefore use a rather rough time-space grid on this interval (i.e. $[t_{N_h-1}, t_{N_h}]$).
Further, only a precise value of  $ \hat v^{m}(t_{N_{h}-1},\cdot)$ will be required for the estimation of $\hat v^{m+1}$ on $[t_{N_{h}-2},t_{N_{h}-1})$ and this is where a fine space grid should be used. Iterating this argument, one can use rather rough time-space grid on each $(t_{i},t_{i+1})$ and concentrate on each $t_{i}$ at which a finer space grid is required. This is the main difference with the usual backward Euler schemes of  \cite{BouchardTouzi,GLW,Zhang} and  the forward Picard schemes of \cite{BenderDenk}.   
\end{remark}

\noindent{\bf Proof of Proposition \ref{prop: erreur discretisation et esperance}.}
	Define 
	\b*
		\tilde v^{m}(\cdot) := (-M_{h_{\circ}})\vee  \Esp{V_{\cdot}(\hat v^m(\tip,\cdot),\hat v^{m-1}) \big| \hat \Fc }\wedge M_{h_{\circ}}.
	\e*
	Then, Lemma \ref{lem:diffprod} below combined with the inequality $|\vp|\le 1$ implies that for all $(t,x)\in [t_i, t_{i+1}) \x \mathbf{X}$,
	\begin{align*}
		&|\tilde v^{m}(t,x)- v^{m}(t,x)|\\
		&\le 
		\E\Big[\sum_{k \in \Kc_{t_{i+1}-t}} 
			\frac{1}{M} V_{t_{i+1}-t}^M  
			\big|\hat v^{m}(\tip,X^{x,(k)}_{\tip})-v^{m}(\tip,X^{x,(k)}_{\tip}) \big|
		\Big| \hat \Fc \Big]
		\\
		&+ 
		\E\Big[ \sum_{k \in \bar \Kc_{t_{i+1}-t} \setminus \Kc_{t_{i+1} -t}}
			2 L_{\varphi} V_{t_{i+1}-t}^M
			\big|  \hat v^{m-1}(T_{k},X^{x,(k)}_{T_{k}})-v^{m-1}(T_{k},X^{x,(k)}_{T_{k}}) \big| 
		\Big| \hat \Fc \Big].
	\end{align*}
	Let us compute the expectation of the first term. Denoting by $\bar \Fc$  the $\sigma$-field generated by the branching processes, we obtain 
	\b*
		&&
		\E \Big[\sum_{k \in \Kc_{t_{i+1}-t}} 
			\frac{1}{M} V_{t_{i+1}-t}^M
			\big|\hat v^{m}(\tip,X^{x,(k)}_{\tip})-v^{m}(\tip,X^{x,(k)}_{\tip}) \big|
		\Big] \\
		&=&
		\E \Big[\sum_{k \in \Kc_{t_{i+1}-t}} 
			\frac{1}{M} V_{t_{i+1}-t}^M
			\E \Big[
			\big|\hat v^{m}(\tip,X^{x,(k)}_{\tip})-v^{m}(\tip,X^{x,(k)}_{\tip}) \big|
		\Big| \bar \Fc \Big] \Big] \\
		&\le&
		\frac{1}{M}   \| \E[| \hat v^m - v^m|] \|_{t_{i+1}} \E \Big[ q_{t_{i+1} -t} V_{t_{i+1} -t}^M \Big]
		~\le~
		\frac{M^1_{h} }{M}   \| \E[| \hat v^m - v^m|] \|_{t_{i+1}}.
	\e*
	Similarly, for the second term, one has
	\b*
		&&
		\E\Big[ \sum_{k \in \bar \Kc_{t_{i+1}-t} \setminus \Kc_{t_{i+1} -t}}
			{ 2} L_{\varphi} V_{t_{i+1}-t}^M
			\big|  \hat v^{m-1}(T_{k},X^{x,(k)}_{T_{k}})-v^{m-1}(T_{k},X^{x,(k)}_{T_{k}}) \big| 
		\Big] \\
		&\le&
		{ 2} L_{\varphi} M^2_{h} \| \E[| \hat v^{m-1} - v^{m-1}|] \|_{t_i}.
	\e*
	Notice that $\|\Esp{|\tilde v^{m} - \hat v^{m}|}\|_{\ti}\le  \Ec(\hat \E)$ by Assumption \ref{assum:hatE}. Hence, 
	\begin{align*}
		\| \E[ | \hat v^{m}- v^{m} |]\|_{\ti}
		&~\le ~
		\Ec(\hat \E) 
		+
		{ 2} L_{\varphi} M^2_{h} \|  \E[ | \hat v^{m-1}- v^{m-1} |] \|_{\ti} \\
		&\;\;\; \;\;\; \;\;\; \;\;\; \;\; +
		\frac{M^1_{h}}{M} \|  \E[ | \hat v^{m}- v^{m} |] \|_{\tip}.
	\end{align*}
	We now  appeal to Proposition \ref{prop: recurrence} to obtain 
	\b*
		\|\E[|\hat v^{m}- v^{m}|] \|_{\ti}
		&\le&
		\Ec(\hat \E) \left(\sum_{i=1}^m C^i +\sum_{i'=2}^{N_{h}-i} \left(\sum_{j_{1}=1}^{m} \cdots \sum_{j_{i'}=1}^{j_{i'-1}} C^{m -j_{i'}}C^{i'-1}\right)\right) \\
		&\le&
		\Ec(\hat \E) (1+N_{h})\frac{(m+N_{h})^{N_{h}}}{N_{h} !} C^{m+N_{h}},
	\e*
	with $C := ({ 2}  L_{\varphi} M^2_{h}) \vee \frac{M^1_{h}}{M} \vee 1$.
\ep

\section{Numerical experiments}
\label{sec:NumericalExamples}
This section is deditacted to some examples ranging from dimension one to five, and showing the efficiency of the methodology exposed above.

In practice, we  modify the algorithm to avoid costly Picard iterations and we propose two versions that  permit to get an accurate estimate 
in only one Picard iteration:\footnote{We omit here the space discretization procedure, for simplicity. It will be explained later on.}
\begin{enumerate}
\item {\bf Method A}:  In the first method, we simply work backward and apply the localization function $\vp_{j}$ to the estimation made on the previous time step. Namely, we replace  \eqref{eq: def hat vnm} by  
\begin{equation}
  \label{eq: def hat vnm_MOD}
  \hat v(t_i,x):=(-M)\vee \Esph{\hat V_{t_i,x}(\hat v(\tip,\cdot),\hat v(\tip,\cdot))}\wedge M.
\end{equation}
Compared to the initial Picard scheme \eqref{eq: def hat vnm}, we expect to need a smaller time  step to reach an equivalent  variance. On the other hand, we do not do any Picard iteration. 
Note that for $i=N_{h}-1$, this corresponds to one Picard iteration with prior given by $g$. 
 For $i=N_{h}-2$, we use the value at $t_{N_{h}-1}$ of the first Picard iteration for the period $[t_{N_{h}-2},t_{N_{h}-1})$ and the initial prior for the last period, etc. 

\begin{remark}\label{rem: picard sur A} 
	This could obviously be complemented by Picard iterations of the form 
\begin{align*}
  \hat v^{m}(t_i,x):=
  (-M)\vee \Esph{\hat V_{t_i,x}(\hat v^{m}(\tip,\cdot),{\rm Lin}[\hat v^{m-1}])}\wedge M,
\end{align*}
with
$$
	 {\rm Lin}[\hat v^{m-1}](t, \cdot) :=\frac{t-\ti}{\tip-\ti}(\hat v^{m-1}(\tip,\cdot)-\hat v^{m-1}(\ti,\cdot))+ \hat v^{m-1}(\ti,\cdot).
$$
In this case, it is not difficult to see that $\hat v^{m}$ coincides with the classical Picard iteration of the previsous sections on $[T_{N_{h}-m},T_{N_{h}}]$   (up to the specific choice of a linear time interpolation). In practice, these additional Picard iterations are not needed, as we will see in the test cases below.
\end{remark}
\item
{\bf Method B}: 
An alternative consists in introducing on each time  discretization mesh  $[t_i, t_{i+1})$ a sub-grid $\hat t_{i,j}= t_i + \hat h j$, $j =0,\ldots,\hat N_h$, with $ \hat t_{i,\hat N_h} =t_{i+1}$
and replace the representation \eqref{eq: practReg} by
\begin{align*}
  \tilde V_{t,x}(\phi, \hat \phi) &:=
                          \Big(\prod_{k \in   \Kc_{\tip-t}}G_{t,x}(\phi,k) \Big) 
                          \Big(\prod_{k \in\bar{\Kc}_{\tip-t}\setminus\Kc_{\tip-t}}A_{t,x}(\hat \phi,k) \Big),\nonumber \\
                         G_{t,x}(\phi,k)&:= \frac{\phi(\tip,X_{\tip-t}^{x,(k)})}{\bar{F}({\tip-t}-T_{k-})}, \nonumber\\
  A_{t,x}(\hat \phi,k)&:= \frac{\sum_{j=1}^{j_{\circ}}a_{j,\xi_k}(X_{T_k}^{x,(k)})\vp_{j}( \hat \phi(\hat  t_{i, \kappa(t+T_{k})},X_{T_k}^{x,(k)})}{p_{\xi_k}\,\rho(\delta_k)}.
\end{align*}
where   $\kappa(t) = \min\{ j >0 : \hat  t_{i,j} \ge t\}$.\\
Then, we evaluate the $v$ function  on  $[t_i, t_{i+1})$ by applying the scheme recursively backward in time for $j=\hat N_h-1, ..,0$:
\begin{equation}
  \label{eq: defresolB} 
  \hat v(\hat t_{i,j},x):=(-M)\vee \Esph{\tilde V_{\hat t_{i,j},x}(\hat v(\tip,\cdot), \hat v)}\wedge M,
\end{equation}
The estimation at date $\hat t_{i,0}= t_i$ is then used as the terminal function  for the previous time step $[t_{i-1}, t_i)$.\\
With this algorithm, we hope that we will be able to take larger time steps $\tip-t_i$ and  to reduce the global computational cost of the global scheme.
Notice that the gain is not obvious: at the level of the inner time steps, the precision does not need to be high as the estimate only serves at selecting the correct local polynomial, however it may be paid in terms of variance. The comment  of Remark \ref{rem: picard sur A} above also applies to this method.
\end{enumerate}

We   first compare the methods A and B on a  simple test case in dimension 1 and then move to more difficult test cases using the most efficient method, which turns
out to be method A.\\
Many parameters   affect the global convergence of the algorithm :
\begin{itemize}
\item The structure of the representation \eqref{eq: f loc poly} used to approximate the driver: we   use quadratic splines or cubic splines \cite{deBoor} to see the
impact of the polynomial approximation. The splines are automatically generated. 
\item The number  $j_{\circ}$ of functions $\vp_j$  used in \eqref{eq: f loc poly}  for the spline representation.
\item The   number of time steps $N_h$ and $\hat N_h$  used  in the algorithm.
\item The grid and the interpolation used on $\Xbf$ for all dates $t_i$, $i=0,\ldots,N_h$. All interpolations are achieved with the StOpt library \cite{GLWa}.
\item The time step   for the Euler scheme   used  to approximate the solution $X$   of \eqref{eq: Diffusion}.
\item The accuracy chosen to estimate the expectations appearing in our algorithm. We compute the empirical standard deviation $\theta$ associated to the Monte Carlo estimation of 
the expectation in \reff{eq: def hat vnm_MOD} or \reff{eq: defresolB}. We try to fix the number $\hat M$ of samples such that   ${\theta}/{\sqrt{\hat M}}$ does not exceed a certain level, fixed below, at each point  of our grid.
\item  The random variables $(\delta_k)_{k}$, which  define the life time of each particle, is chosen to follow an exponential distribution with parameter $0.4$ for all of the studies below.
\end{itemize}

\subsection{A first simple test case in dimension one}
In this section, we compare the methods A and B on the following  toy example. 
Let us set $\Xbf:=[\underline x,\bar x]$ with $\underline x=\pi/8$ and $\bar x=7\pi/8$, and consider the solution $X$ of \reff{eq: Diffusion} with 
$$
\mu(x)=0.1\x (\frac{\pi}{2}-x)\;\mbox{ and }\; \sigma(x):=0.2\x(\bar x-x)(x-\underline x).
$$
	
We then consider the driver 
\begin{flalign*}
f(t,x,y)= & \mu(x) \hat f(y) + (-\alpha +\frac12\sigma(x)^2) y
\end{flalign*}
where 
\begin{flalign*}
\hat f(y) = & \left(\sqrt{e^{2 \alpha (T-t)} -y^2}\1_{|y|\le  \bar y e^{ \alpha (T-t)}}+ \tilde f(y) \1_{|y|>  e^{ \alpha (T-t)} \bar y}\right) \\
\tilde f(y) = & \sqrt{e^{2 \alpha (T-t)}-\bar y^2 e^{ 2\alpha (T-t)}} -  \frac{\bar y e^{\alpha (T-t)}  }{\sqrt{e^{2 \alpha (T-t)}-\bar y^2 e^{ 2\alpha (T-t)}}} (|y| -\bar y e^{ \alpha (T-t)})
\end{flalign*}
with $\bar y:=\cos(\underline x)$. The solution is $v(t,x)= e^{-\alpha (T-t)} \cos(x)$ for $g(x)= \cos(x)$. 

Notice the following points :
\begin{itemize}
\item  this case favors the branching method because  $g$  and  $v$ are bounded by one, meaning that the variance due to a product of pay-off is bounded by one too.
\item  in fact the domain of interest is  such that $|y| \le \bar y e^{ \alpha (T-t)}$ but we need to have a smooth approximation of the driver at  all the  dates on the
domain $[-\bar y, \bar y]$.
\end{itemize}
We   take the following parameters: $T=1$, $\alpha =0.5$, the time step used for the interpolation scheme is $0.4$,  and we use the modified monotonic quadratic interpolator defined in \cite{warin}. The Euler scheme's discretization step  used  to simulate the Euler scheme of  $X$  is  equal to $0.002$, the number of simulations $\hat M$ used  for the branching method is taken such that $\theta/\sqrt{\hat M} < 0.002$ and limited to $10.000$ (recall that $\theta$ is the empirical standard deviation). At last the truncation parameter $M$ is set to one. \\
\vs1

\begin{figure}[H]
\centering
\includegraphics[height=6cm,width=10cm]{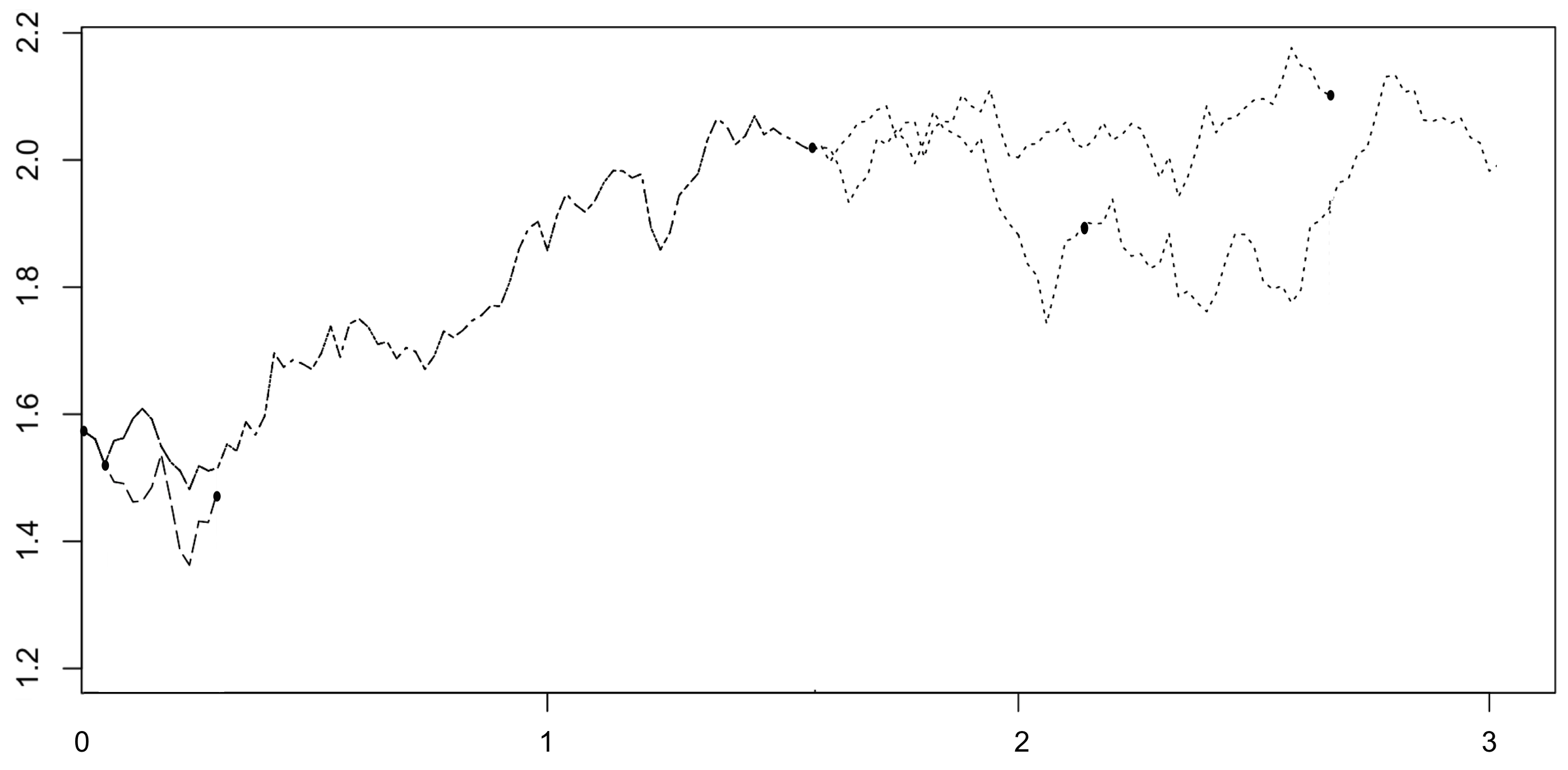}
\caption{A typical simulated path of the branching diffusion starting from $\pi/2$ on $[0,3]$. Bullets denote  branching or killing times. }
\label{fig: path ex}
\end{figure}

A typical path of the branching diffusion  is provided in Figure \ref{fig: path ex}. It starts from $x = \pi/2$ at $t=0$. 

\begin{figure}[H]
\begin{minipage}[b]{0.5\linewidth}
  \centering
 \includegraphics[width=0.9\textwidth]{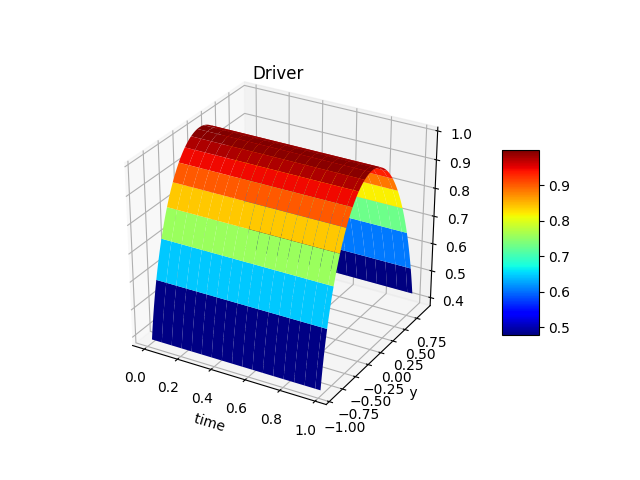}
 \caption*{The driver}
 \end{minipage}
\begin{minipage}[b]{0.5\linewidth}
  \centering
 \includegraphics[width=0.9\textwidth]{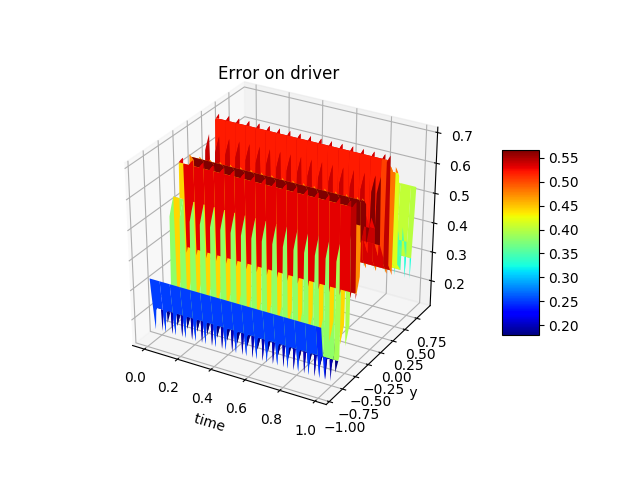}
 \caption*{Error on the driver due to the spline representation.}
 \end{minipage}
\caption{The driver $\hat f$ and its quadratic  spline representation error for 20 splines. } 
\label{fig:driverBruno}
\end{figure}

To estimate the driver $\hat f$, we use a quadratic spline: on Figure \ref{fig:driverBruno} we plot  $\hat f$  on $ [0,1] \times [-\bar y, \bar y]$  and the error obtained with a $20$ splines representation. 
Notice that this driver has a high Lipschitz constant around $-\bar y$ and $\bar y$. 
\vs1

\begin{figure}[H]
\centering
\includegraphics[height=4cm,width=8cm]{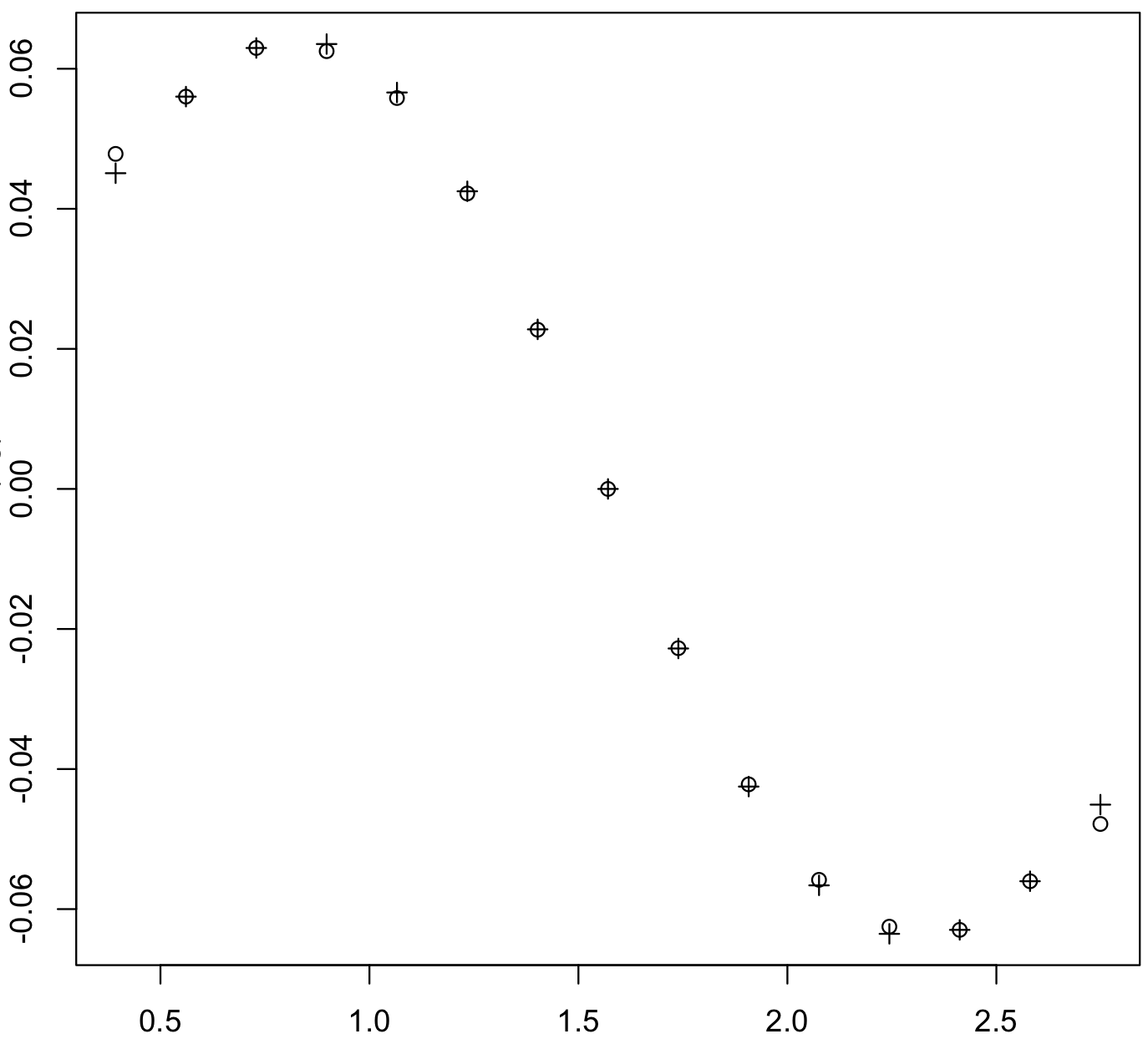}
\caption{Approximation of the driver at $t=T$ with only three local polynomials -  Crosses: $f(\cdot,\cos)$. Circles:  $f_{\ell_{\circ}}(\cdot,\cos,\cos)$.}
\label{fig: f approx}
\end{figure}

As already mentionned, we do not try to optimize the local polynomial representation, but instead generate the splines automatically. Our motivation is that the method should work in an industrial context, in which case a hand-made approximation might be complex to construct.  Note however that one could  indeed, in this test case, already achieve a very good precision with only three local polynomials as shown in Figure  \ref{fig: f approx}. Recalling Remark \ref{rem: localisation}, this particular approximation would certainly be more efficient, in particular if the probability of reaching the boundary points, at which the precision is not fully satisfactory, is small.
\vs1

\begin{figure}[H]
\begin{minipage}[b]{0.5\linewidth}
  \centering
 \includegraphics[width=0.9\textwidth]{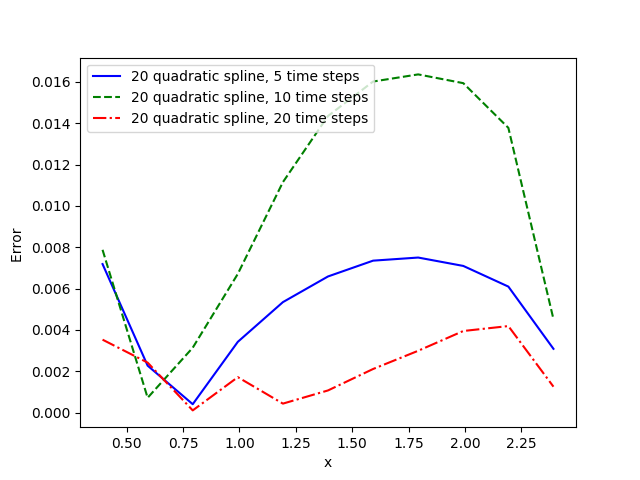}
 \caption*{\small Error in absolute value for different\\ time discretization.}
 \end{minipage}
\begin{minipage}[b]{0.5\linewidth}
  \centering
 \includegraphics[width=0.9\textwidth]{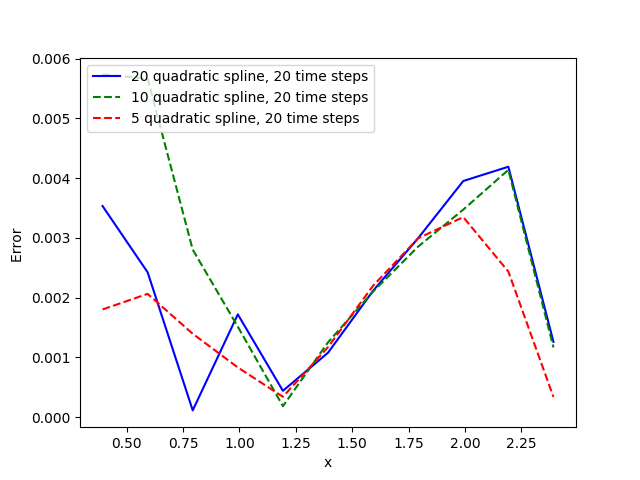}
 \caption*{\small Error in absolute value for different\\ number of splines.}
 \end{minipage}
\caption{Results for the toy case using method A } 
\label{fig::methodA}
\end{figure}

\begin{figure}[H]
  \centering
 \includegraphics[width=0.5\textwidth]{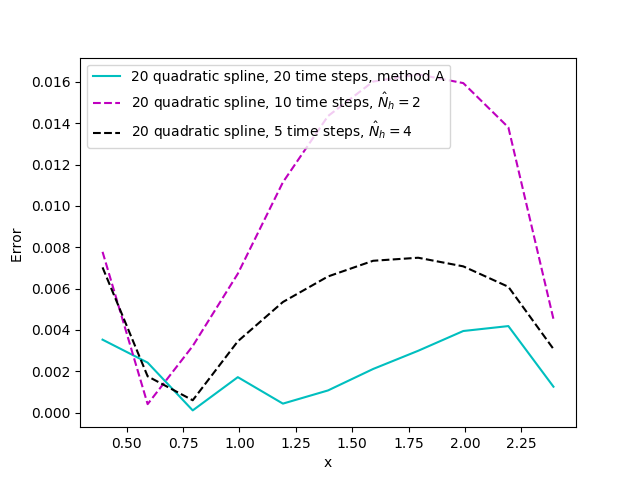}
 \caption{\label{fig::methodB}Error in absolute value for method B for different $N_h$ and $\hat N_h$, keeping  $N_h \hat N_h=20$.}
\end{figure}

In Figure \ref{fig::methodA}, we give the results obtained by method A for different values of $N_h$, the number of time steps,   and of  the number of splines used to construct $\hat f$. As shown on the graph, the error with $N_h=20$ and $20$ splines is below $0.004$, and even with $5$ time steps the results are very accurate. Besides the results obtained are very stable with the number of splines used. 
On Figure \ref{fig::methodB}, we give the results obtained with the method B for different values of $N_h$ and $\hat N_h$ keeping the number  $N_h \hat N_h$ constant,
equal to $20$.
Comparing   Figure \ref{fig::methodA} and Figure \ref{fig::methodB}, we see that the results obtained with the two methods are similar, for a same total number of time steps. This shows that method B is  less efficient, as it gives similar results but at the cost of additional  computations.
This is   easily explained by the fact that the variance increases a lot with the size of the time step,  so as to become the first source of error.
In the sequel of the article, we will therefore concentrate on method A.\\
The time spent for the best results, obtained with method A,  $20$ time steps and $20$ splines,  is less than $1.05$  seconds on a regular (old indeed) laptop.

\subsection{Some more difficult  examples}
We show in this section that the  method works well even in more difficult situations, in particular when the boundary condition  $g$   is not bounded by one. This will be at the
price of a higher variance, that is compensated by an increase of the computational cost.\\
We now take  $\Xbf:=[0,2]^d$,  with   
$$
\mu(x)= U \x (\1 -x)\;\mbox{ and }\; \sigma(x):= V \prod_{i=1}^d (2-x_i) x_i {\rm I_{d}},
$$
where $V=0.2$, $U=0.1$ and ${\rm I_{d}}$ is the identity matrix.\\

We will describe the drivers later on. Let us just immediatly mention that  we use the modified monotonic quadratic interpolator defined in \cite{warin} for tensorized interpolations in dimension $d$. When using sparse grids, we apply the quadratic interpolator of \cite{bungartz3}. The spatial discretization  used for all tests is defined with
\begin{itemize}
\item a step 0.2 for the full grid interpolator,
\item a level 4  for the sparse grid interpolator \cite{bungartz1}.
\end{itemize} 
We take a very thin time step of $0.000125$ for the  Euler scheme  of $X$.
The number of simulations $\hat M$  is chosen such ${\theta}/{\sqrt{\hat M}} < 0.00025$ and limited to $200.000$,  so that the error
reached can be far higher than $0.00025$ for high time steps. This is due to fact that the empirical standard deviation $\theta$  is large for some points near the boundary of $\Xbf$. We finally use a truncation parameter   $M=50$, it  does not appear to be very relevant numerically.

\subsubsection{ A first one dimensional example}
\begin{figure}[H]
\begin{minipage}[b]{0.5\linewidth}
  \centering
 \includegraphics[width=0.9\textwidth]{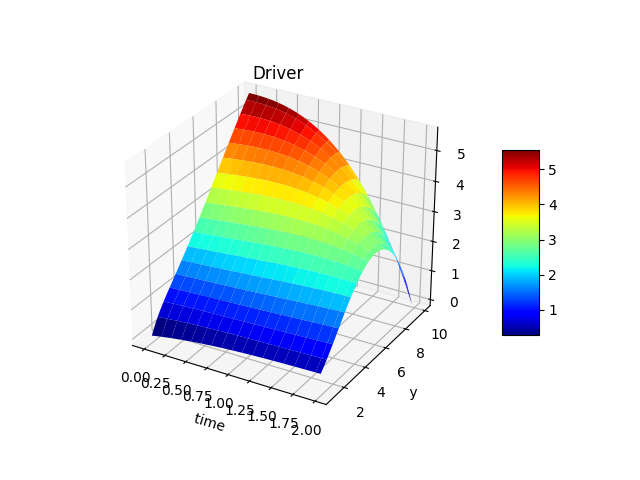}
 \caption*{The driver}
 \end{minipage}
\begin{minipage}[b]{0.5\linewidth}
  \centering
 \includegraphics[width=0.9\textwidth]{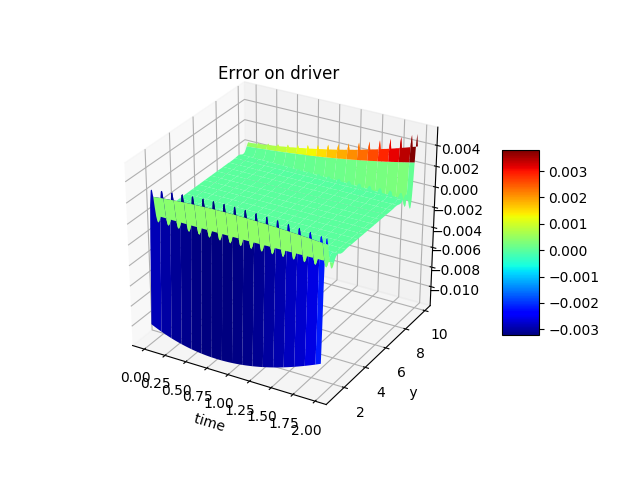}
 \caption*{Error on the driver due to the spline representation.}
 \end{minipage}
\caption{The driver and its cubic spline representation error with 10 splines for the first difficult case. } 
\label{fig:driver1}
\end{figure}

\begin{figure}[!h]
\begin{minipage}[b]{0.5\linewidth}
  \centering
  \includegraphics[width=\textwidth]{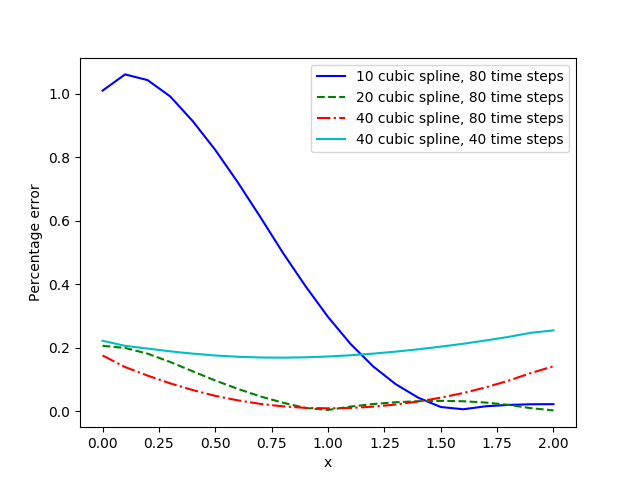}
  \caption{ Cubic spline method.}
  \label{fig:case1Cubic}
\end{minipage}
\begin{minipage}[b]{0.5\linewidth}
  \centering
  \includegraphics[width=\textwidth]{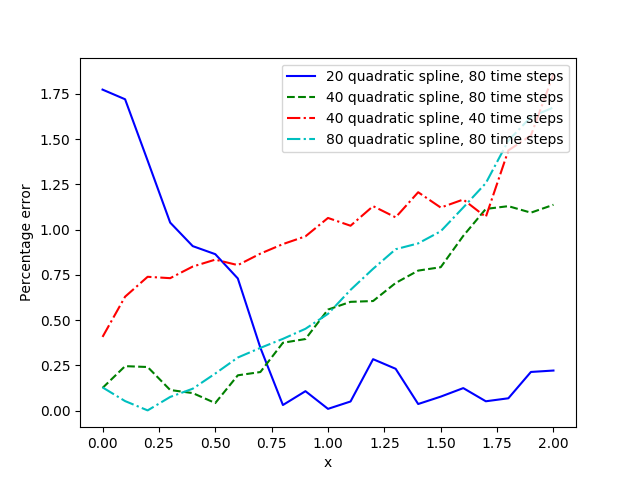}
  \caption{ Quadratic spline method.}
  \label{fig:case1Quad}
\end{minipage}
\caption*{Percentage error on the approximation of  the solution $v(0,.)$ for different time discretization and spline discretization. }
\end{figure}

In this part,  we use the following time dependent driver  for a first one dimensional example:
\begin{flalign}
f(t,y)= & y (\frac{1}{2}  - \frac{V^2}{2C^2} [ \phi(t,T,y) (2C-\phi(t,T,y))]^2 - U(C -\phi(t,T,y))), \nonumber
\end{flalign}
with $\phi(t,T,y)= \log(y)-\frac{T-t}{2}.$\\
We use a time discretization of $1.000$ time steps to represent the time dependency of the driver.
On Figure \ref{fig:driver1},  we provide the driver and the cubic spline error associated to 10 splines.
Figure  \ref{fig:case1Cubic} corresponds to  cubic splines (with an approximation on each mesh with a polynomial of degree 3),
while    Figure \ref{fig:case1Quad} corresponds to  quadratic  splines  (with an approximation on each mesh with a polynomial of degree 2).\\

On this example, the cubic spline approximation appears to be far more efficient than the quadratic spline.
In order to get a very accurate solution when using sparse grids, with a maximum  error below $0.2\%$, it is necessary to have a high number of splines (at least $20$) and a high number of time steps, meaning  that the  high variance of the method for the highest time step has prevented the algorithm to converge with the maximum number of samples imposed.

\subsubsection{A second one dimensional example}
\begin{figure}[H]
\begin{minipage}[b]{0.5\linewidth}
  \centering
 \includegraphics[width=0.9\textwidth]{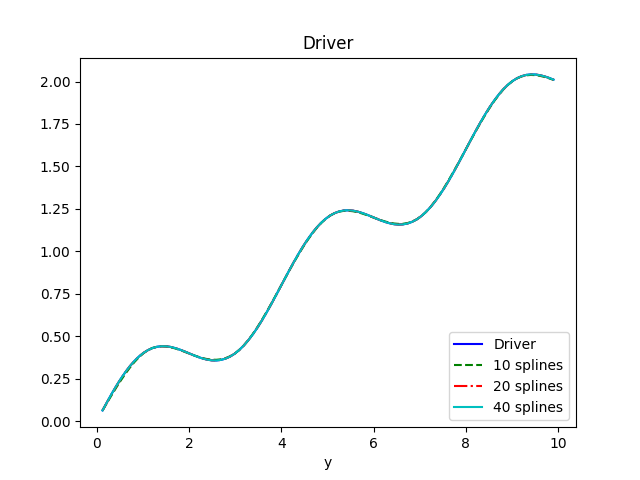}
 \caption*{The driver}
 \end{minipage}
\begin{minipage}[b]{0.5\linewidth}
  \centering
 \includegraphics[width=0.9\textwidth]{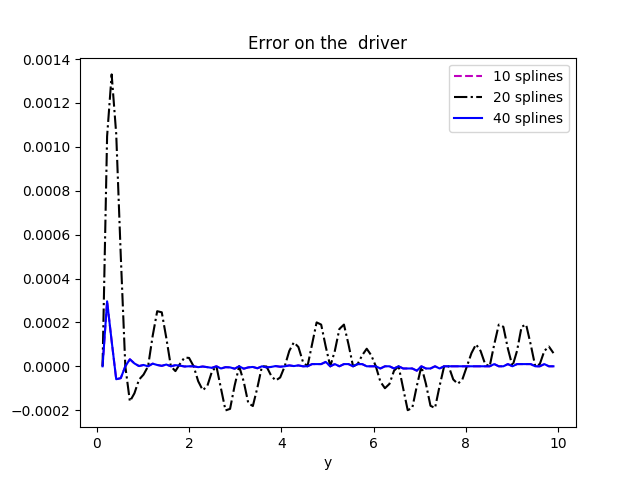}
 \caption*{Error on the driver due to the spline representation.}
 \end{minipage}
\caption{The driver and its cubic spline representation. } 
\label{fig:driver2}
\end{figure}

\begin{figure}[H]
\begin{minipage}[b]{0.5\linewidth}
  \centering
  \includegraphics[width=\textwidth]{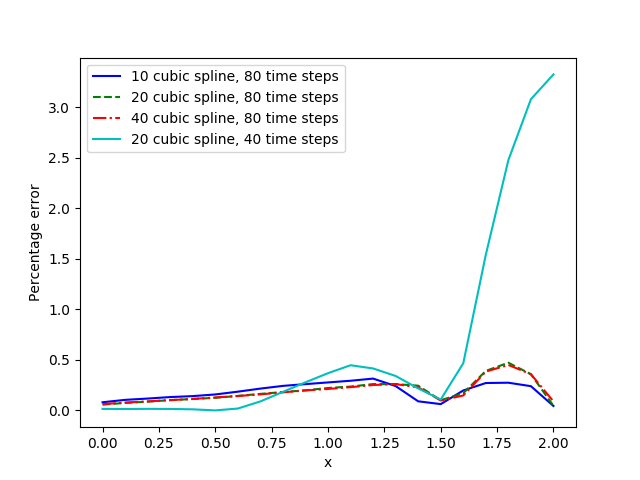}
  \caption{ Cubic spline method.}
  \label{fig:case2Cubic}
\end{minipage}
\begin{minipage}[b]{0.5\linewidth}
  \centering
  \includegraphics[width=\textwidth]{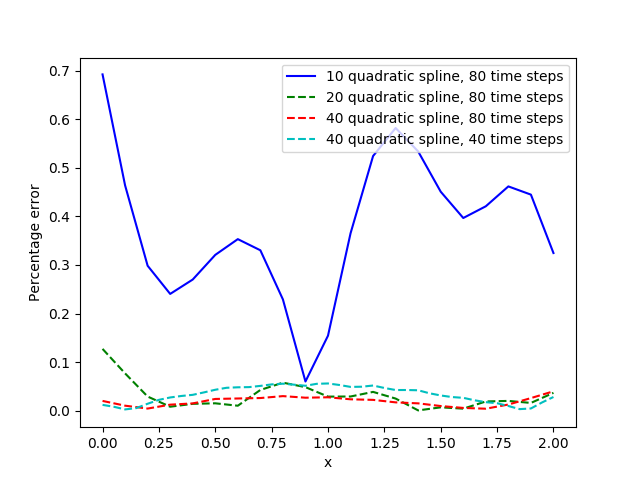}
  \caption{ Quadratic spline method.}
  \label{fig:case2Quad}
\end{minipage}
\caption*{Percentage error on the aproximation of the solution $v(0,.)$ for different time discretizations and numbers of splines. }
\end{figure}

We now consider  the  driver 
\begin{flalign}
\label{secondDriver}
f(x,y)= f_1(y)+ f_2(x),
\end{flalign}
with
\begin{flalign}
f_1(y) =& \frac{2}{10} (y + \sin(\frac{\pi}{2} y)),\label{f1f2}  \\
f_2(x) =&  \frac{1}{2}-(\frac{2}{10} + C \mu(x)) - \frac{\sigma(x)^2 c^2}{2} e^{C x +\frac{T-t}{2}} - \frac{2}{10} \sin(\frac{\pi}{2} e^{c x  +\frac{T-t}{2}}).\nonumber
\end{flalign}
Figure \ref{fig:driver2} shows $f_1$  and the cubic spline error associated to different numbers of  splines.
On Figure \ref{fig:case2Cubic} and  \ref{fig:case2Quad}, we give, in percentage,  the error obtained when using cubic and quadratic splines for different discretizations.
Globally,  the quadratic interpolator appears to provide better results for both coarse and thin  discretizations.
With $40$ time steps, the cubic approximation generate errors up to $3\%$ at the boundary point $x=2$: reviewing the results at each time step, we checked that the convergence of the Monte Carlo is not achieved near the boundary $x=2$, with Monte Carlo errors up to $0.03$ at each step.
\vs1

Finally, note that the convergence of the method is related to the value of  the quadratic and cubic coefficients of the spline representation, that we want to be as small as possible.

\subsubsection{Multidimensional results}
\begin{figure}[H]
\begin{minipage}[b]{0.5\linewidth}
  \centering
 \includegraphics[width=\linewidth] {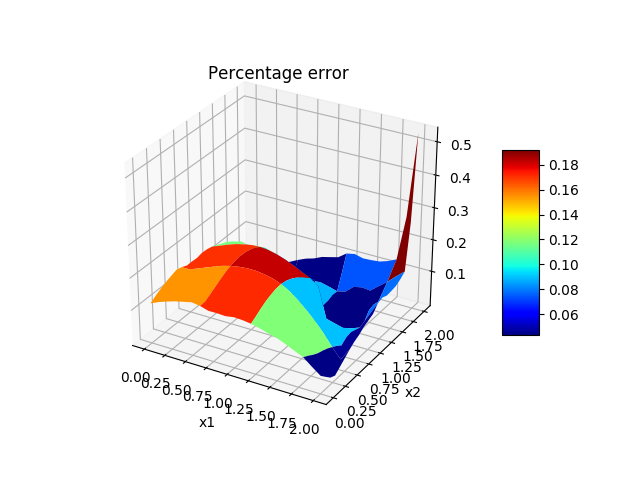}
 \caption*{10 splines}
 \end{minipage}
\begin{minipage}[b]{0.5\linewidth}
  \centering
 \includegraphics[width=\textwidth]{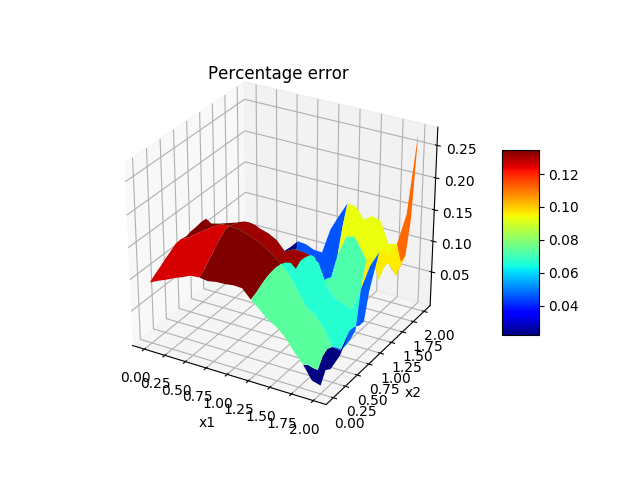}
 \caption*{20 splines}
 \end{minipage}
\begin{minipage}[b]{0.5\linewidth}
  \centering
 \includegraphics[width=\textwidth]{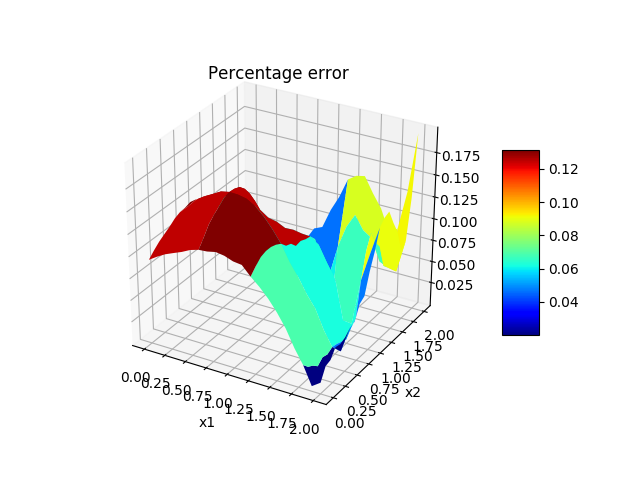}
 \caption*{40 splines}
 \end{minipage}
\caption{Error in dimension 2 for $80$ time steps with cubic splines. } 
\label{fig:secondTest2D}
\end{figure}

\begin{figure}[H]
\begin{minipage}[b]{0.5\linewidth}
  \centering
  \includegraphics[width=\linewidth] {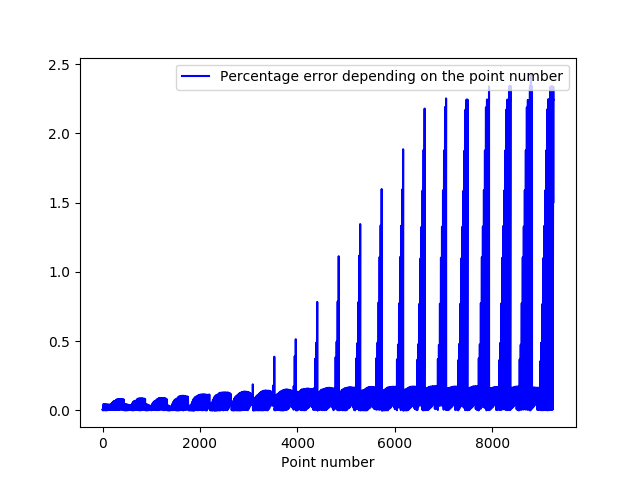}
 \caption*{40 splines, 80 time steps.}
 \end{minipage}
\begin{minipage}[b]{0.5\linewidth}
  \centering
 \includegraphics[width=\textwidth]{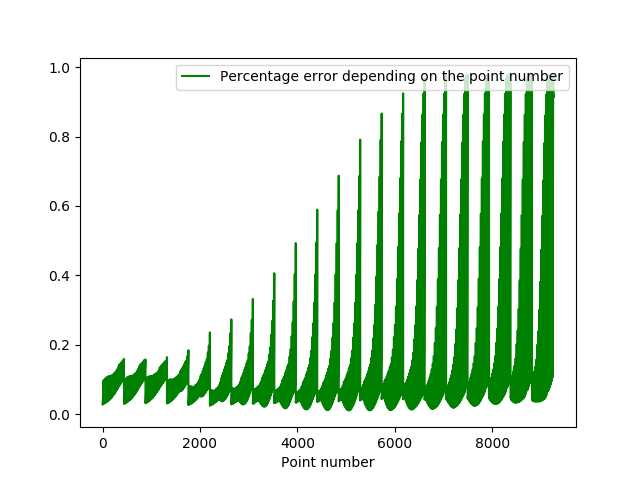}
 \caption*{80 splines, 80 time steps}
 \end{minipage}
\begin{minipage}[b]{0.5\linewidth}
  \centering
 \includegraphics[width=\textwidth]{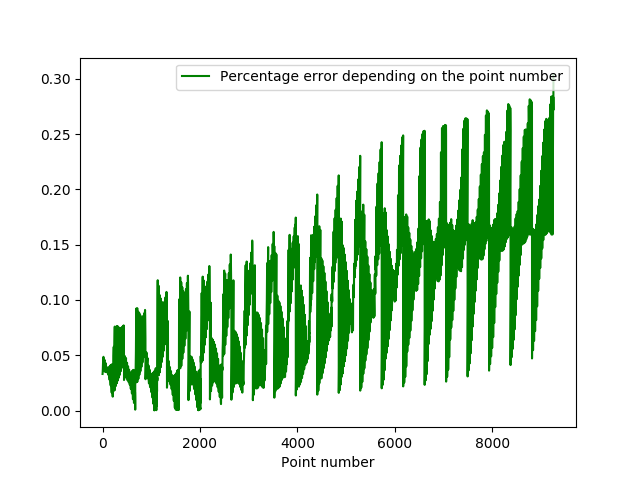}
 \caption*{80 splines, 120 time steps}
 \end{minipage}
\begin{minipage}[b]{0.5\linewidth}
  \centering
 \includegraphics[width=\textwidth]{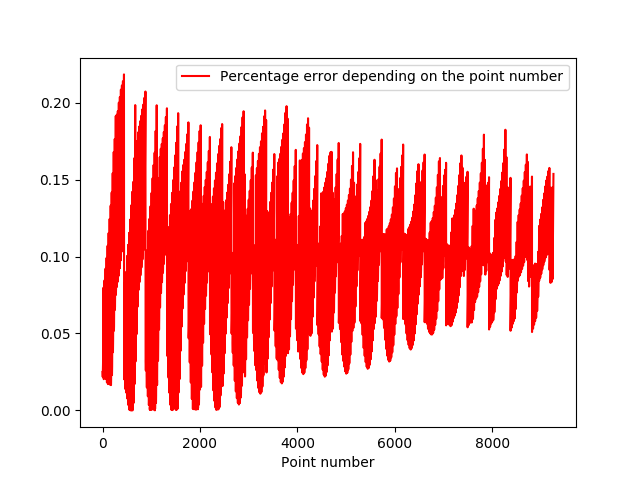}
 \caption*{80 splines, 160 time steps}
 \end{minipage}
\caption{Error in dimension 3 for different time steps  and spline numbers, with cubic splines. } 
\label{fig:secondTest3D}
\end{figure}

In this section, we keep the driver  in the form \eqref{secondDriver} with   $f_1$ as    in \eqref{f1f2}, but we now generalize the definition of $f_{2}$:
\begin{flalign*}
 f_2(x)=& \frac{1}{2}-(\frac{2}{10} + \frac{C}{d}  \sum_{i=1}^d x_i ) - \frac{\sigma^{11}(x)^2 c^2}{2d} e^{\frac{C}{d}  \sum_{i=1}^d x_i +\frac{T-t}{2}} \\
 &-   \frac{2}{10} \sin(\frac{\pi}{2} e^{\frac{C}{d}  \sum_{i=1}^d x_i +\frac{T-t}{2}}) 
\end{flalign*}
In this section, we only consider  cubic splines.
\paragraph{Results with full grids}
Figure \ref{fig:secondTest2D} describes the results   in dimension 2 for $80$ time steps and different numbers of splines. Once again the number of splines used is relevant for the convergence.
Figure \ref{fig:secondTest3D} corresponds to dimension $3$ for different numbers of splines and different time discretizations: the error is plotted as a function of the point number using a classical Cartesian numeration.
Once again, the results are clearly improved when  we increase the number of splines and increase the number of time steps.  
 
\begin{figure}[H]
\centering
 \includegraphics[width=0.5\textwidth]{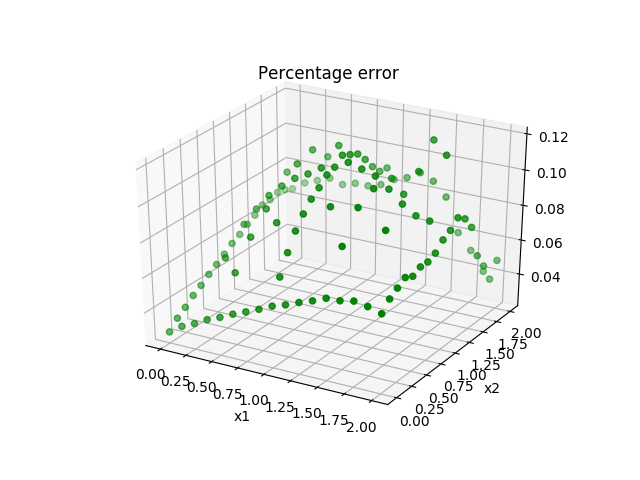}
\caption{Error of the quadratic sparse grid of level 4  with 80 time steps}
\label{fig:sparse2D}
\end{figure}

\begin{figure}[H]
\begin{minipage}[b]{0.5\linewidth}
  \centering
  \includegraphics[width=0.9\linewidth] {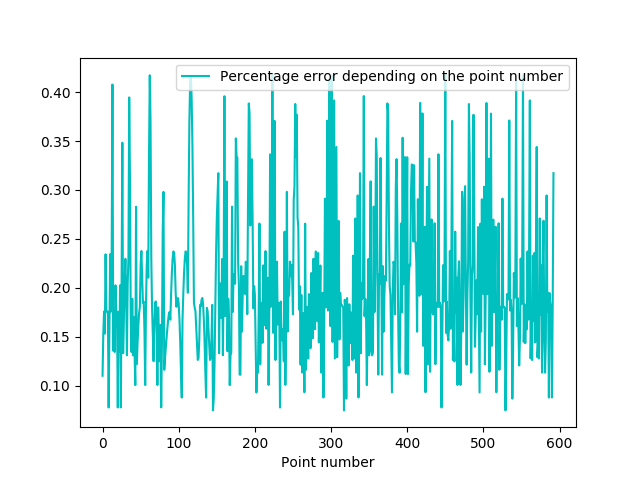}
 \caption*{3D, 80 splines, 160 time steps.}
 \end{minipage}
\begin{minipage}[b]{0.5\linewidth}
  \centering
 \includegraphics[width=0.9\textwidth]{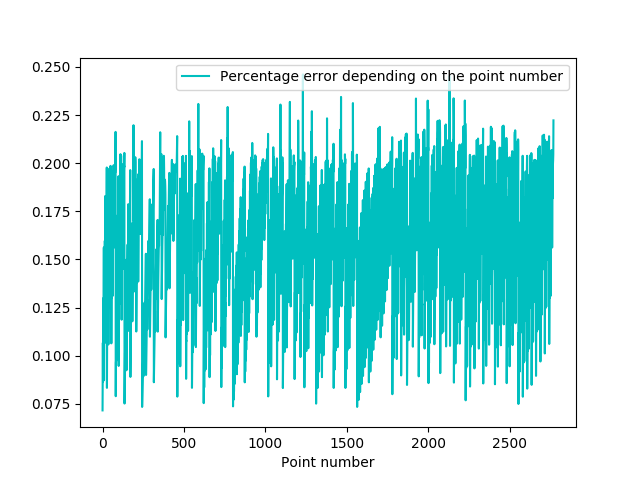}
 \caption*{4D, 80 splines, 120 time steps.}
 \end{minipage}
\begin{minipage}[b]{0.5\linewidth}
  \centering
 \includegraphics[width=0.9\textwidth]{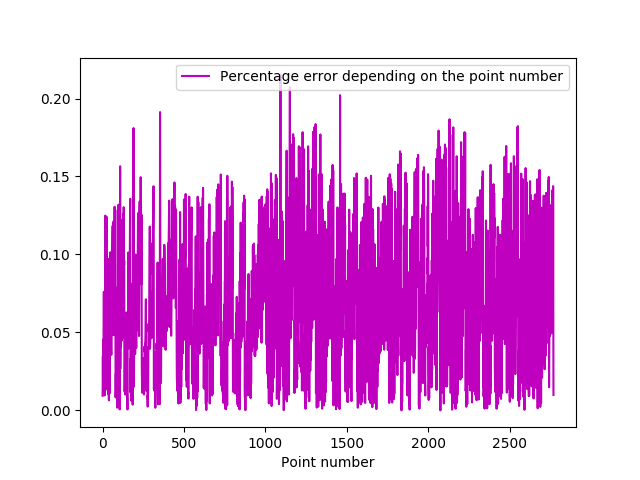}
 \caption*{4D, 80 splines, 160 time steps.}
 \end{minipage}
\begin{minipage}[b]{0.5\linewidth}
  \centering
 \includegraphics[width=0.9\textwidth]{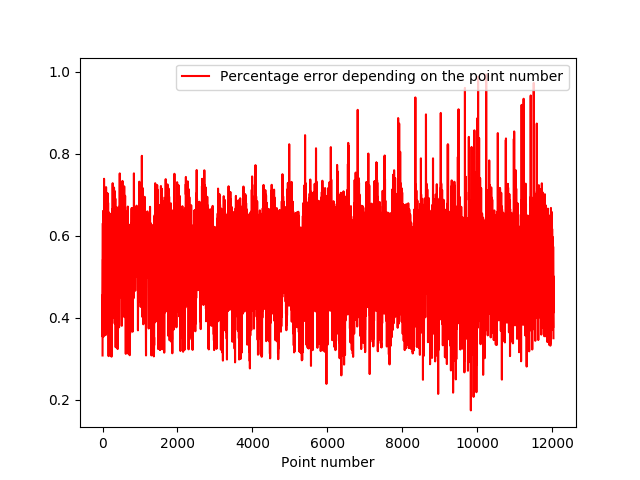}
 \caption*{5D , 80 splines, 160 time steps.}
 \end{minipage}
\caption{Error in dimensions 3, 4 and 5  with a cubic spline approximation and a quadratic sparse interpolator. } 
\label{fig:secondTestSparse}
\end{figure}

\paragraph{Towards higher dimension}
As the dimension increases, the algorithm    is subject to the curse of dimensionality. This is due to the $d$-dimensional interpolation: the number
of points used is $n^d$ if $n$ is the number of points in one dimension.
One way to surround this  is to use sparse grids \cite{bungartz1}. The sparse grid methodology permits  to get an interpolation error nearly as good as with full grids for a cost increasing slowly with the dimension, whenever the solution is
smooth enough. It is based on 
some special interpolation points. According to  \cite{bungartz1,bungartz3}, if the function $w$  to be interpolated is null at the boundary and admits derivatives such that 
$\sup_{\alpha_i \in \{2,3\}} \left \{ || \frac{\partial^{\alpha_1+..+\alpha_d} \;w}{\partial x_1^{\alpha_1} ... \partial x_d^{\alpha_d}}||_{\infty} \right \} < \infty$,
then the interpolation error due to the quadratic sparse interpolator $I^2$ is given 
\begin{flalign}
\label{interpSparse}
|| w - I^2(w)||_\infty =  O(n^{-3} log(n)^{d-1}).
\end{flalign}
An effective sparse grids implementation  is given in \cite{GLWa} and more details on sparse grids   can be found in \cite{GLWa2}.\\
On Figure \ref{fig:sparse2D}, we plot the error obtained with a $2$-dimensional  sparse grid, for $80$ time steps. 
In dimension  3, 4, and 5,   the error obtained with the spline of level 4 is given in Figure \ref{fig:secondTestSparse}. 
Once again, we are able to be very accurate.
\begin{remark}
On Figure \ref{fig:secondTestSparse}, the error is plotted as a function of the point number. Due to the special structure of the sparse grid  points numeration, no special pattern appears (as on Figure \ref{fig:secondTest3D}) but the maximum error is still located at the boundary of $\Xbf$.
\end{remark}
From a practical point of view, the algorithm can be easy parallelized: at each time step each, points can be affected to one processor.
It can therefore be speeded-up (linearly) with respect to the number of cores used. 

\begin{figure}[H]
 \centering
 \includegraphics[width=0.5\textwidth]{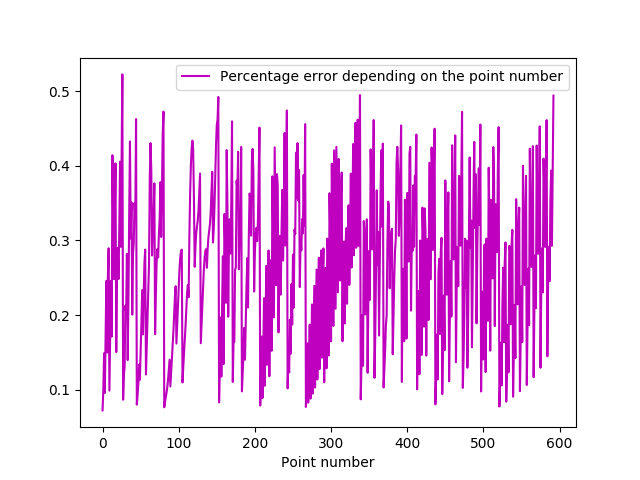}
\caption{Error in dimension 3, 80 splines, 160 time steps and less tight  parameters.}
\label{fig:sparse3DRela}
\end{figure}

All results have been obtained using a cluster with a MPI implementation for the parallelization.
The time spend using the sparse grid interpolator in 3D, a Monte Carlo accuracy   fixed to  $0.0005$  and for an Euler scheme with time step  $0.001$, is $2.800$ seconds on a cluster with  8 processors  using 112 cores. The error curve  can be found in Figure \ref{fig:sparse3DRela}.

\begin{remark}
When the solution has not the required regularity to reach the interpolation error \eqref{interpSparse}, 
one can implement a local adaptation of the grid using  a classical estimation of the local error (based on the hierarchical surplus) \cite{griebel,bungartz1}. It is also possible to use the dimension adaptive method   \cite{gerstner}, which aims at refining a whole dimension when a higher interpolation error in this  dimension is detected. 
\end{remark}
\begin{remark}
The methodology developed here is very similar in spirit to the Semi-Lagrangian method used in \cite{warin,warin1} using full   or    sparse grids:
the deterministic scheme starting from one point of the grid is replaced by a Monte Carlo one. Classically, the error of such a scheme is decomposed into a time discretization error and an interpolation error. This is the same in our scheme but,  to the time discretization error due to the use of the scheme  \eqref{eq: def hat vnm_MOD}, is added a Monte Carlo error associated to the branching scheme.
\end{remark}

\newpage
\appendix
\section{Appendix}

\subsection{Technical lemmas}

\begin{lemma} \label{lemm:estim_ODE}
	The ordinary differential equation $\eta'(t) = \sum_{\ell = 0}^{\ell_{\circ}} 2 C_{\ell_{\circ}} \eta(t)^{\ell}$
	with initial condition $\eta(0) = M > 0$
	has a unique solution on   $[0, h_{\circ}]$ for 
	\be \label{eq:def_delta}
		h_{\circ} := \frac{(\ell_{\circ} - 1) (1-M)_+ + (1 \vee M)^{-(\ell_{\circ}-1)}}{(\ell_{\circ} +1)(\ell_{\circ} -1) 2C_{\ell_{\circ}}}.
	\ee
	Moreover, it is bounded on $[0, h_{\circ}]$ by
	\begin{equation} \label{eq:def_M_delta}
		M_{h_\circ} 
		:= 
		\max \!\Big(1, 
			\big( (1 \vee M)^{1 - \ell_{\circ}} \!+\! (\ell_{\circ} \!-\!1) (1\!-\! M)^+ \!-\! h_{\circ} \ell_{\circ} (\ell_{\circ}\!-\!1) 2C_{\ell_{\circ}} \big)^{(1-\ell_{\circ})^{-1}}
		 \Big).
	\end{equation}
	Consequently, one has, for all $t \in [0, h_{\circ}]$,
	\be \label{eq:branch_bounded_by_Mh}
		\E\Big[
		\Big( \prod_{k \in \Kc_t} \frac{M}{\bar F(t - T_{k-})} \Big)
		\Big( \prod_{k \in \bar \Kc_t \setminus \Kc_t} \frac{ 2C_{\ell_{\circ}}}{p_{\xi_k} \rho(\delta_k)} \Big)
		\Big] 
		~\le~
		\eta(t) ~\le~ M_{h_\circ}.
	\ee
\end{lemma}
\proof
	$\mathrm{i)}$ We first claim that
	\be \label{eq:claimMh}
		\int_M^{M_{h_\circ}} \frac{dy}{2C_{\ell_{\circ}}(1+y+ \cdots + y^{\ell_{\circ}})}
		\ge h_{\circ}.
	\ee
	Then, for every $t \in [0, h_{\circ}]$, there is some constant $M(t) \le M_{h_\circ} < \infty$ such that
	$$
		\int_M^{M(t)} \frac{dy}{2C_{\ell_{\circ}}(1+y+ \cdots + y^{\ell_{\circ}})}  ~=~ t = \int_0^t ds.
	$$
	This means that $(M(t))_{t \in [0, h_{\circ}]}$ is a bounded solution (and hence the unique solution) of  
	$\eta'(t) = \sum_{\ell = 0}^{\ell_{\circ}} 2 C_{\ell_{\circ}} \eta(t)^{\ell}$
	with initial condition $\eta(0) = M > 0$.
	In particular, it is bounded by $M_{h_\circ}$.

	\vspace{1mm}

	$\mathrm{ii)}$ Let us now prove \eqref{eq:claimMh}.
	Notice that $y^k \le 1 \vee y^{\ell_{\circ}}$ for any $y \ge 0$ and $k = 0, \cdots, \ell_{\circ}$. Then, it is enough to prove that
	\be \label{eq:ODE_interm_estim}
		\int_M^{M_{h_{\circ}}} \Big(1 \wedge \frac{1}{y^{\ell_{\circ}}} \Big) dy \ge h_{\circ} ( \ell_{\circ} +1) 2C_{\ell_{\circ}}.
	\ee
	By direct computation, the l.h.s.~of \eqref{eq:ODE_interm_estim} equals  
	$$
		(M_{h_\circ} - M) \1_{\{M_{h_\circ} \le 1\}}
		+
		\Big(
			(1-M)^+ + \frac{1}{\ell_{\circ}-1} \big( (1\vee M)^{1-\ell_{\circ}} - M_{h_\circ}^{1-\ell_{\circ}} \big)
		\Big) \1_{\{M_{h_\circ} > 1\}}.
	$$
	When $h_{\circ}$ satisfies \eqref{eq:def_delta}, it is easy to check that \eqref{eq:ODE_interm_estim} holds true.
	
	\vspace{1mm}
	
	$\mathrm{iii)}$ We now prove \eqref{eq:branch_bounded_by_Mh}.
	Recall that  $\bar \Kc_t^n$ denotes  the collection of all particles in   $\bar \Kc_t$ of generation $n$. Set
	$$
		\chi^n_t := \Big( \prod_{k \in \cup_{j=1}^n \Kc^j_t} \frac{M}{\bar F(t - T_{k-})} \Big)
		\Big( \prod_{k \in \cup_{j=1}^n (\bar \Kc^j_t \setminus \Kc^j_t)} \frac{ 2C_{\ell_{\circ}}}{p_{\xi_k} \rho(\delta_k)} \Big)
		\Big( \prod_{k \in \bar \Kc^{n+1}_t} \eta(t - T_{k-}) \Big).
	$$
	Since $\bar \Kc^n_t$ has only finite number of particles, the random variable $\chi^n_t$ is uniformly bounded.
	Then by exactly the same arguments as in \eqref{eq:branch_rep_app1} and \eqref{eq:branch_rep_app2} below,
	and by repeating this argument over $n$, one has
	\b*
		\eta(t) 
		=
		M + \int_0^t \sum_{\ell=0}^{\ell_{\circ}} 2 C_{\ell_{\circ}} \eta(s)^{\ell} ds
		=
		\E\big[ \chi^1_t \big] 
		=
		\E \big[ \chi^n_t \big],
		~~~\forall n \ge 1.
	\e*
	It follows by Fatou Lemma that
	$$
		\E\Big[
		\Big( \prod_{k \in \Kc_t} \frac{M}{\bar F(t - T_{k-})} \Big)
		\Big( \prod_{k \in \bar \Kc_t \setminus \Kc_t} \frac{ 2C_{\ell_{\circ}}}{p_{\xi_k} \rho(\delta_k)} \Big)
		\Big]
		=
		\E \big[ \lim_{n \to \infty} \chi^n_t \big]
		~\le~
		\lim_{n \to \infty} \E[ \chi^n_t] = \eta(t).
	$$
\qed

\vspace{2mm}

For completeness, we provide here the proof the representation formula of Proposition  \ref{propeq: representation Ynm} and of the technical lemma that was used in the proof of Proposition \ref{prop: erreur discretisation et esperance}.

\begin{proposition}\label{prop: representation Ynm} The representation formula of Proposition  \ref{propeq: representation Ynm} holds. 
\end{proposition}

\proof We only provide the proof on $[t_{N_{h}-1},T]$, the general result is obtained by induction. It is true by construction when   $m$ is equal to $0$. Let us now fix $m\ge 1$. 

\vspace{1mm}

First, Lemma \ref{lemm:estim_ODE} shows that the random variable $V^m_{t,x}$ is integrable.

\vspace{1mm}

Next, Set $(1)+:=\{(1,j), j \le \ell_{\circ}\}\cap \bar \Kc_{T}$ and define $\Kc_t(1):=\Kc_t \cap  (1)+$ and $\bar \Kc_t(1):=\bar \Kc_t \cap (1)+$. 
	{For ease of notations, we write $X^x := X^{x, ((1))}$.}
	Then, for all $(t,x) \in [t_{N_h-1}, T] \x \R^d$, 
\begin{align*}
\E[V^{m}_{t,x}]&=\E\left[\frac{g(X_{T-t}^{x})}{\bar{F}(T-t)}\1_{\{T_{(1)}\ge T-t\}}\right]\\
&+\E\left[\1_{\{T_{(1)}< T-t\}}\frac{\sum_{j=1}^{j_{\circ}}a_{j,\xi_{(1)}}(X_{T_{(1)}}^{x})\vp_{j}(v^{m-1}(t+T_{(1)},X_{T_{(1)}}^{x}))}{p_{\xi_{(1)}}\,\rho(\delta_{(1)})}R^{m}_{t,x}\right]
\end{align*}
where 
$$
R^{m}_{t,x}:=
\Big(\prod_{k \in   \Kc_{T-t}(1)}G_{t,x}(k) \Big) \Big( \prod_{k \in\bar{\Kc}_{T-t}(1)\setminus\Kc_{T-t}(1)}A^{m}_{t,x}(k) \Big)
$$
satisfies 
$$
\E[R^{m}_{t,x}|\Fc_{T_{(1)}}]=\prod_{k\in (1)+} v^{m} \big( {t+}T_{(1)},X^{t,x}_{T_{(1)}} \big)= \big[v^{m}( {t+}T_{(1)},X^{x}_{T_{(1)}}) \big]^{\xi_{(1)}},
$$
by \reff{hyp eq: independance}. On the other hand, \reff{eq: def bar F} and \reff{hyp eq: independance} imply 
\begin{align} \label{eq:branch_rep_app1}
\E\left[\frac{g(X_{T-t}^{x})}{\bar{F}(T-t)}\1_{\{T_{(1)}\ge T-t\}}\right]&=\E[g(X_{T-t}^{x})]
\end{align}
and 
\begin{align} \label{eq:branch_rep_app2}
&\E\left[ \1_{\{T_{(1)}<T-t\}} \frac{\sum_{j=1}^{j_{\circ}}a_{j,\xi_{(1)}}(X_{T_{(1)}}^{x})\vp_{j}( v^{m-1}(t+T_{(1)},X_{T_{(1)}}^{x}))}{p_{\xi_{(1)}}\,\rho(\delta_{(1)})}[v^{m}({t+}T_{(1)},X^{x}_{T_{(1)}})]^{\xi_{(1)}}\right] \nonumber \\
&=\E\left[  \int_{0}^{T-t} \frac{\sum_{j=1}^{j_{\circ}}a_{j,\xi_{(1)}}(X^{x}_{s})\vp_{j}(v^{m-1}(t+s,X_{s}^{x}))}{p_{\xi_{(1)}}} [v^{m}( {t+}s,X^{x}_{s})]^{\xi_{(1)}} ds \right] \nonumber \\
&=\E\left[   \int_{0}^{T-t} \sum_{j=1}^{j_{\circ}}\sum_{\ell \le \ell_{\circ}} a_{j,\ell}(X_{{s}}^{x})\vp_{j}(v^{m-1}(t+s,X_{{s}}^{x}))  [v^{m}( {t+}s,X^{x}_{s})]^{\ell }ds \right] \nonumber \\
&=\E\left[   \int_{0}^{T-t} f_{\ell_{\circ}}(X_{{s}}^{x},v^{m}(t+s,X_{{s}}^{x}),v^{m-1}(t+s,X_{{s}}^{x}))ds \right].
\end{align}
Combining the above implies that 
$$
v^{m}(t,X_{t})=\E\left[g(X_{T})+\int_{t}^{T} f_{\ell_{\circ}}(X_{s},v^{m}(s,X_{s}),v^{m-1}(s,X_{s}))ds \Big| \Fc_{t}\right],
$$
and the required result follows by induction. 
\ep 
\vs2
 
\begin{lemma}\label{lem:diffprod}
Let $(x^{i},y^{i})_{i \le I}$ be a sequence of real numbers. 
Then, 
	\begin{align*}
		\left|\prod_{i=1}^{I} x^{i}-\prod_{i=1}^{I} y^{i}\right|
		\le
		\sum_{i \in I}  
		\Big( 
			|x^i - y^i| \prod_{j \neq i}  \max(|x^j|, |y^j|)
		\Big).
	\end{align*}
\end{lemma}
\proof It suffices to observe that 
	\begin{align*}
		\prod_{i=1}^{I} x^{i}-\prod_{i=1}^{I} y^{i}
		~=~
		(x^{1}-y^{1})\prod_{i=2}^{I} x^i
		+ y^1 \Big( \prod_{i=2}^I x^i - \prod_{i=2}^{I} y^{i} \Big),
	\end{align*}
	and to proceed by induction.
	\ep

 \begin{proposition}\label{prop: recurrence} 
 	Let $c_{1},c_{2},c_{3}\ge 0$, and let $(u^{i}_{m})_{m\ge 0,i\ge 0}$ be a sequence such that 
 	$$
 		u^{i}_{m}\le c_{1}u^{i}_{m-1}+c_{2}u^{i+1}_{m}+c_{3} \;\mbox{ for } m\ge 1, i<N_{h}.
 	$$
 	Then
	\begin{align*}
		u^{i}_{m}
		\le& 
		c_{1}^{m}u^{i}_{0}+\sum_{i'=1}^{N_{h}-i} \left(\sum_{j_{1}=1}^{m}\sum_{j_{2}=1}^{j_{1}} \cdots \sum_{j_{i'}=1}^{j_{i'-1}} c_{1}^{m}c_{2}^{i'}
		u^{i+i'}_{0}\right)\\
		&+
		c_{3}  \left( \sum_{i=1}^m c^i_{1} + \sum_{i'=2}^{N_{h}-i} \Big(\sum_{j_{1}=1}^{m}\sum_{j_{2}=1}^{j_{1}} \cdots \sum_{j_{i'}=1}^{j_{i'-1}} c_{1}^{m-j_{i'}}c_{2}^{i'-1} \Big) \right).
	\end{align*}
 \end{proposition}
\proof We have 
\begin{align*}
u^{i}_{m}\le& (c_{1})^{m}u^{i}_{0}+\sum_{j=1}^{m} (c_{1})^{m-j} (c_{2}u^{i+1}_m +c_{3} ).
\end{align*}
The required result then follows from a simple induction.\ep

\subsection{More on the error analysis for the abstract numerical approximation}

	The regression error  $\Ec(\hat \E)$ in Assumption \ref{assum:hatE}
	  depends essentially on the regularity of $v^m$. 
	Here we prove that $v^m(t,x)$ is H\"older in $t$ and Lipschitz in $x$ under additional conditions,
	and provide some estimates on the corresponding coefficients.
	Given $\phi: [0,T] \x \R^d \to \R$, denote
	\b*
		[\phi]_{\ti} ~:=~ 
		\sup_{(t,x)\ne (t', x') \in [\ti,\tip] \x \Xbf } 
		\frac{|\phi(t,x)-\phi(t',x')|}{|t-t'|^{\frac12}+|x-x'|}.
	\e*
	Since $(\mu,\sigma)$ is assumed to be Lipschitz, it is clear that there exists $L_{X}>0$ such that 
	for all $(t,x),(t',x')\in [0,T]\x \Xbf$,
	\be\label{eq: estimate holder X}
		\|X^{x}_{t}-X^{x'}_{t'}\|_{\Lb^{2}}\le L_{X}  \Big( \sqrt{|t' - t|} + |x'-x| \Big).
	\ee
  
\begin{proposition}\label{prop: holder vm} 
	Suppose that $x \mapsto g(x)$ and $x \mapsto f_{\ell_{\circ}}(x,y,y')$ are uniformly Lipschitz with Lipschitz constants $L_g$ and $L_f$.
	Let $\beta$ and $\lambda_1, \lambda_2 > 0$ such that
	$\frac{L_2}{\lambda^2_2} T < 1$ and $\beta \ge 2L_1 + L_f \lambda_1^2 + L_2 \lambda^2_2$,
	then for all $m \ge 1$ and $ i \le N_{h}$,
	\b*
		[v^{m}]_{\ti} ~\le~ L_{v} 
		&:=&
		(1 + L_X) L_X \sqrt{\big(L_g^2 + \frac{L_f}{ \beta \lambda^2_1} \big)   e^{\beta T}/ \big(1 - \frac{L_2}{\lambda_2^2}T \big)}  \\
		&&+~
		2 (1+ \ell_{\circ}) C_{\ell} ( 1 \vee (M_{h_\circ})^{\ell_{\circ}}) \sqrt{h_{\circ}}.
	\e*
\end{proposition}  
\proof
	For ease of notations, we provide the proof for $t=0$ only.  
	
	$\mathrm{i)}$ Let $x_1, x_2 \in \R^d$ and $Y^{m,1} := v^m(\cdot, X^{x_1})$, $Y^{m,2} := v^m(\cdot, X^{x_2})$,
	and denote $\Delta Y^m := Y^{m,1}  - Y^{m,2}$, $\Delta X := X^{x_1} - X^{x_2}$,
	where $X^{x_1}$ (resp. $X^{x_2}$) denotes the solution of SDE \eqref{eq: Diffusion} with initial condition $X_0 = x_1$ (resp. $X_0 = x_2$).
	Using the same arguments as in the proof of Theorem \ref{thm: main},
	it follows that, for any $\beta \ge 2L_1 + L_f \lambda_1^2 + L_2 \lambda^2_2$, one has
	\be \label{eq:estim_Lip_interm}
		\E [e^{\beta t} (\Delta Y_t^{m+1})^2] 
		&\le& 
		\E [e^{\beta T} (\Delta Y_T^{m+1})^2] 
		 ~+~ \frac{L_f}{\lambda_1^2} \E \Big[ \int_t^T e^{\beta s} |\Delta X_s |^2 ds \Big]
		\nonumber \\
		&&
		+ \frac{L_2}{\lambda^2_2}\E \Big[ \int_t^T e^{\beta s} (\Delta Y^{m}_s)^2 ds \Big]
	\ee
	and then
	\b*
		\E \Big[ \int_0^T e^{\beta t} (\Delta Y^{m+1}_t)^2 dt \Big]
		&\le&
		T \E [e^{\beta T} (\Delta Y_T^{m+1})^2]
		+
		T \frac{L_f}{\lambda_1^2} \E \Big[ \int_0^T e^{\beta s} |\Delta X_s|^2 ds \Big] \\
		&&
		+~ T \frac{L_2}{\lambda_2^2}
		\E \Big[ \int_0^T e^{\beta t} (\Delta Y^m_t)^2 dt \Big]\\
		&\le&
		T e^{\beta T} \Big( L_{g}^{2} + \frac{L_f}{\beta \lambda_1^2} \Big) L_{X}^{2} |x_{1}-x_{2}|^{2} \\
		&&
		+~
		T\frac{L_2}{\lambda_2^2} 
		\E \Big[ \int_0^T e^{\beta t} (\Delta Y^m_t)^2 dt \Big].
	\e*
	Since $\frac{L_2}{\lambda^2_2} T < 1$, this induces that
	$$
		\E \Big[ \int_0^T e^{\beta t} (\Delta Y^{m+1}_t)^2 dt \Big]
		~\le~
		\frac{T e^{\beta T} \big(L_g^2 + \frac{L_f}{ \beta \lambda^2_1} \big)  L_X^2 |x_1 - x_2|^2}{1 - \frac{L_2}{\lambda^2_2} T}.
	$$
	Plugging the above estimates into \eqref{eq:estim_Lip_interm}, it follows that
	$$
		(\Delta Y_0^m )^2
		~\le~
		\hat L_v^2 |x_1 - x_2|^2,
		~~~\mbox{with}~~
		\hat L_v^2 :=\frac{\big(L_g^2 + \frac{L_f}{ \beta \lambda^2_1} \big) L_X^2 e^{\beta T}}{1 - \frac{L_2}{\lambda^2_2} T}.
	$$

	$\mathrm{ii)}$ For the H\"older property of $v^m$, it is enough to notice that for $t \le h_{\circ}$,
	\b*
		|v^m(0,x) - v^m(t,x)| 
		&\le& 
		\E \Big[ |v^m(t, X^x_t) - v^m(t, x)| + \int_0^t |f(X^x_s, Y^m_s, Y^{m-1}_s)| ds \Big] \\
		&\le&
		\hat L_v L_X \sqrt{ t} + 2 (1+ \ell_{\circ}) C_{\ell} ( 1 \vee (M_{h_\circ})^{\ell_{\circ}}) t,
	\e*
	where the last inequality follows from the Lipschitz property of $v^m$ in $x$ and the fact that $Y^m$ is uniformly bounded by $M_{h_\circ}$.
	We hence conclude the proof.
\ep

\section*{Ackowledgements}

This work has benefited from the financial support of the Initiative de Recherche ``M\'ethodes non-lin\'eaires pour la gestion des risques financiers'' sponsored by AXA Research Fund.

\vspace{1mm}

Bruno Bouchard and Xavier Warin   acknowledges the financial support of ANR project CAESARS (ANR-15-CE05-0024).

\vspace{1mm}

Xiaolu Tan  acknowledges the financial support of the ERC 321111 Rofirm, the ANR Isotace, and the Chairs Financial Risks (Risk Foundation, sponsored by Soci\'et\'e G\'en\'erale) and Finance and Sustainable Development (IEF sponsored by EDF and CA).

\end{document}